\documentclass{amsart}

\usepackage{amsthm,amsfonts,amsmath,amssymb,latexsym,epsfig}
\usepackage{upref,eucal}

\def \cH{{\mathcal H}}
\def \h{\hat{\ }}

\def \bZ{{\mathbb Z}}
\def \bG{{\mathbb G}}

\def \<{\langle}
\def \>{\rangle}

\def \cF{\mathcal F}

\def \cX{\mathcal X}

\def \cF{\mathcal F}

\def \d{\delta}
\def \bZ{{\mathbb Z}}
\def \bF{{\mathbb F}}
\def \bC{{\mathbb C}}
\def \bR{{\mathbb R}}
\def \bA{{\mathbb A}}
\def \cO{\mathcal O}
\def \ra{\rightarrow}
\def \bF{{\mathbb F}}
\def \bQ{{\mathbb Q}}
\def \bQ{{\mathbb Q}}

\def \bP{{\mathbb P}}
\def \bG{{\mathbb G}}
\def \bW{{\mathbb W}}

\newtheorem{theorem}{Theorem}[section]

\theoremstyle{definition}

\theoremstyle{remark}

 \numberwithin{equation}{section}
\theoremstyle{problem}
\newtheorem{problem}[theorem]{Problem}
\newtheorem{example}[theorem]{Example}
\begin{document}
\title{Differential calculus with integers}
\author{Alexandru Buium}

\maketitle

\medskip

\centerline{\it University of New Mexico, Albuquerque, USA}

\begin{abstract}
Ordinary differential equations have an arithmetic analogue in which functions are replaced by numbers and the derivation operator is replaced by a Fermat quotient operator.
In this survey we explain  the main motivations, constructions, results, applications, and open problems of the theory.
\end{abstract}

\bigskip

\bigskip

The main purpose of these notes is to show how one can develop
 an arithmetic analogue of differential calculus in which differentiable functions $x(t)$ are replaced by integer numbers  $n$ and the derivation operator $x\mapsto \frac{dx}{dt}$ is  replaced by the Fermat quotient operator $n \mapsto \frac{n-n^p}{p}$, where $p$ is a prime integer. The Lie-Cartan geometric theory of differential equations  (in which solutions are smooth maps) is then replaced by a  theory  of ``arithmetic differential equations" (in which solutions are integral points of algebraic varieties). 
 In particular the differential invariants of groups in the Lie-Cartan theory are replaced by ``arithmetic differential invariants"
 of correspondences  between algebraic varieties. 
 A number of  applications to diophantine geometry  over number fields and to classical modular forms will be explained. 

This program was initiated in \cite{char} and pursued, in particular,  in  \cite{pjets}-\cite{adel3}; for an exposition  of some of these  ideas we refer to the monograph \cite{book}. 
We shall restrict ourselves here to the {\it ordinary differential} case.
 For the {\it partial differential} case we refer to \cite{laplace,pde,pdemod,forms}.
Throughout these notes we assume familiarity with the basic concepts of  algebraic geometry and differential geometry; some of the standard material is being reviewed, however, for the sake of introducing notation, and ``setting the stage".
The notes are organized as follows. The first section presents some classical background, the main concepts of the theory, a discussion of the main motivations, and a comparison with other theories. The second section presents a sample of the main results. The third section presents a list of open problems.

\medskip

{\bf Acknowledgement}.  The author is indebted to HIM for support during part of the semester on Algebra and Geometry in Spring 2013. These notes are partially based on lectures given at the IHES in Fall 2011 and MPI in Summer 2012 when the author was partially supported by IHES and MPI respectively.
 Partial support was also received from  the NSF through grant DMS 0852591.

\section{Main concepts}

 \subsection{Classical analogies}
 The analogies between functions and numbers have played a key role in the
development of  modern number theory. Here are some classical analogies.
All facts in this subsection are well known and entirely classical; we review them only in order to introduce notation and put things in perspective.

 \subsubsection{Polynomial functions}  The 
 ring $\bC[t]$ of polynomial functions with complex coefficients is analogous to
the ring $\bZ$ of integers. The field of rational functions $\bC(t)$ is then analogous to the field of rational numbers $\bQ$. In $\bC[t]$ any non-constant polynomial is a product of linear  factors. In $\bZ$ any integer different from $0,\pm 1$ is up to a sign a product of prime numbers. To summarize
$$\bC \subset \bC[t]\subset \bC(t)$$
are analogous to 
$$\{0,\pm 1\}\subset \bZ\subset \bQ$$

\subsubsection{Regular functions}  More generally rings $\cO(T)$ of regular functions on complex algebraic affine non-singular curves $T$ are analogous to rings of integers $\cO_F$ in number fields $F$. 
Hence curves $T$ themselves are analogous to schemes $Spec\ \cO_F$. Compactifications
$$T \subset \overline{T}=T\cup \{\infty_1,...,\infty_n\}\stackrel{\ }{\simeq} \text{(compact Riemann surface of genus $g$)}$$
are analogous to ``compactifications"
$$Spec\ \cO_F\subset \overline{Spec\ \cO_F}=(Spec\ \cO_F) \cup \frac{Hom(F,\bC)}{\text{conjugation}}$$

\subsubsection{Formal functions} The inclusions
$$\bC \subset \bC[[t]] \subset \bC((t))$$
(where $\bC[[t]]$ is  the ring of power series and $\bC((t))$ is the ring of Laurent series) are analogous to the inclusion
$$\{0\}\cup \mu_{p-1}=\{c\in \bZ_p;c^p=c\}\subset \bZ_p\subset \bQ_p$$
(where $\bZ_p$ is ring of  $p$-adic integers and $\bQ_p=\bZ_p[1/p]$). Recall that 
$$\bZ_p=\lim_{\leftarrow} \bZ/p^n\bZ= \{\sum_{n=0}^{\infty} c_i p^i; c_i\in \{0\}\cup\mu_{p-1}\}$$
So $\{0\}\cup \mu_{p-1}$ plays the role of ``constants" in $\bZ_p$. Sometimes we need more ``constants" and we are led to consider, instead, the inclusions:
$$\{0\}\cup \bigcup_{\nu} \mu_{p^{\nu}-1} \subset \widehat{\bZ_p^{ur}}\subset \widehat{\bZ_p^{ur}}[1/p]$$
where 
$$\widehat{\bZ_p^{ur}}=\bZ_p[\zeta; \zeta^{p^{\nu}-1}=1, \nu\geq 1]\h=
\{\sum_{i=0}^{\infty} c_i p^i; c_i\in \{0\}\cup \bigcup_{\nu} \mu_{p^{\nu}-1}\}.$$
Here the upper hat on a ring $A$ means its $p$-adic completion:
$$\widehat{A}:=
\lim_{\leftarrow} A/p^nA.$$
So in the latter case the monoid $\{0\}\cup \bigcup_{\nu} \mu_{p^{\nu}-1}$
should be viewed as the set of ``constants" of  $\widehat{\bZ_p^{ur}}$; this is consistent with the ``philosophy of the field with one element" to which we are going to allude later.
Let us say that a ring is  a {\it local $p$-ring} if it is a discrete valuation ring  with maximal ideal generated by a prime $p\in \bZ$. Then
 $\bZ_p$ and  $\widehat{\bZ_p^{ur}}$ are local $p$-rings.
 Also for any local $p$-ring $R$ we denote by $k=R/pR$ the residue field and by $K=R[1/p]$ the fraction field of $R$.
 Sometimes we will view local $p$-rings as analogues of rings $\bC\{x\}$ of germs of analytic functions on Riemann surfaces and even as analogues of rings of global analytic (respectively $C^{\infty}$ functions) on a Riemann surface $T$ (respectively on a $1$-dimensional real manifold $T$, i.e. on a circle $S^1$ or $\bR$).

\subsubsection{Topology}    Fundamental groups 
of complex curves (more precisely Deck transformation groups of normal covers $T'\ra T$ of Riemann surfaces) have, as  analogues, Galois groups $G(F'/F)$ of
normal extensions number fields $F\subset F'$. The genus of a Riemann surface has an analogue for number fields defined in terms of ramification.
All of this is very classical. There are other, less classical, topological analogies like the one  between primes in $\bZ$ and nodes in $3$-dimensional real manifolds \cite{knots}. 

\subsubsection{Divisors} The group of divisors 
$$Div(\overline{T})=\{\sum_{P\in \overline{T}} n_P P; n_P\in \bZ\}$$
on a non-singular complex algebraic curve $\overline{T}$ is analogous to the group of divisors
$$Div(\overline{Spec\ \cO_F})=\{\sum \nu_P P; \nu_P\in \bZ \ \text{if $P$ is finite}, \nu_P\in \bR\ \text{if $P$ is infinite}\}$$
One can attach divisors to rational functions $f$ on $\overline{T}$ ($Div(f)$ is the sum of poles minus the sum of zeroes); similarly one can attach divisors to elements  $f\in F$. In both cases one is lead to a ``control" of the spaces of $f$s that have a
 ``controlled" divisor (the Riemann-Roch theorem).
One also defines in both settings divisor class groups. In the geometric setting the divisor class group of $\overline{T}$ is an  extension of $\bZ$ by the Jacobian
 $$Jac(\overline{T})=\bC^g/(\text{period lattice of $\overline{T}$})$$ where $g$ is the genus of $\overline{T}$. In the number theoretic setting divisor class groups can be interpreted as ``Arakelov class groups"; one recaptures, in particular, the usual class groups $Cl(F)$. Exploring usual class groups ``in the limit", when one adjoins roots of unity, leads to Iwasawa theory. We will encounter Jacobians later in relation, for instance, to the Manin-Mumford conjecture. This conjecture (proved by Raynaud) says that if one views $\overline{T}$ as embedded into $Jac(\overline{T})$ (via the ``Abel-Jacobi map") the the intersection of $\overline{T}$ with the torsion group of $Jac(\overline{T})$ is a finite set. 
This particular conjecture does not seem to have an analogue for numbers.

 \subsubsection{Families} Maps $$M \ra T$$ of complex algebraic varieties, complex analytic,  or real smooth manifolds, where $\dim T=1$, are analogous to arithmetic schemes i.e. schemes of finite type $$X \ra Spec\ R$$ where $R$ is either  the ring of integers $\cO_F$ in a number field $F$ or a complete local $p$-ring respectively.
 Note however that in this analogy one ``goes arithmetic only half way": indeed for $X \ra Spec\ R$ the basis is arithmetic yet the fibers are still geometric.  One can attempt to ``go arithmetic all the way" and find an  analogue of $M\ra T$ for which both the base and the fiber are ``arithmetic"; in particular one would like to have an analogue of $T \times T$ which is ``arithmetic" in two directions. This is one of the main motivations in the search for $\bF_1$, the ``field with one element".

 \subsubsection{Sections}
 The set of sections 
 $$\Gamma(M/T)=\{s:T\ra M; \pi \circ \sigma=1\}$$ of a map $\pi:M \ra T$ is analogous to the set 
 $$X(R)=\{s:Spec\ R\ra X; \pi\circ s=1\}$$ of $R$-points of $X$ where $\pi:X\ra Spec\ R$ is the structure morphism. This analogy suggests that finiteness conjectures for sets of the form $X(R)$, which one makes in Diophantine geometry, should have as analogues finiteness conjectures for sets of sections $\Gamma(M/T)$. A typical example of this phenomenon is the Mordell conjecture (Faltings' theorem) saying that if $X$ is an algebraic curve of ``genus" $\geq 2$, defined by polynomials with coefficients in a number field $F$, then the set  $X(F)$ is finite. Before the proof of this conjecture Manin \cite{manin} proved a parallel finiteness result for $\Gamma(M/T)$ where $M\ra T$ is a ``non-isotrivial" morphism from an algebraic surface to a curve, whose fibers have genus $\geq 2$. Manin's proof involved the consideration of differential equations with respect to vector fields on $T$. Faltings' proof went along completely different lines. This raised the question whether one can develop a theory of differential equations in which one can differentiate numbers.

\medskip

All  these examples of analogies are classical; cf. work of Dedekind, Hilbert, Hensel, Artin, Weil, Lang, Tate, Iwasawa, Grothendieck, and many others.

\subsection{Analogies proposed in \cite{char}-\cite{adel3}}
 One thing that  seems to be missing from the classical
picture is a counterpart, in number theory,
 of the differential calculus (in particular of differential equations) for functions.
Morally the question is  whether one can meaningfully consider (and successfully use) ``arithmetic differential equations" satisfied by
numbers.
In our research on the subject \cite{char} - \cite{adel3} we proposed  such a theory based on the following sequence of analogies:

\subsubsection{Derivatives} The derivative operator $\d_t=\frac{d}{dt}:\bC[t]\ra \bC[t]$ is analogous to the {\it Fermat quotient operator}
$\d=\d_p:\bZ \ra \bZ$, $\d_p a=\frac{a-a^p}{p}$.
 More generally the derivative operator $\d_t=\frac{d}{dt}:C^{\infty}(\bR)\ra C^{\infty}(\bR)$ (with $t$ the coordinate on $\bR$) is analogous to the operator $\d=\d_p:R \ra R$ on a complete local $p$-ring $R$, $\d_p \alpha=\frac{\phi(\alpha)-\alpha^p}{p}$ (where $\phi:R \ra R$ is a fixed homomorphism lifting the $p$-power Frobenius map on $R/pR$). The map $\d=\d_p$ above is, of course, not a derivation but, rather, it satisfies the following conditions:
 $$
 \begin{array}{lcl}
 \d(1) & = & 0\\
 \d(a+b) &= &\d(a)+\d(b)+C_p(a,b)\\
\d(ab)&=&a^p\d(b)+b^p\d(a)+p\d(a)\d(b),\end{array}$$
where $C_p(x,y)\in \bZ[x,y]$ is the polynomial $C_p(x,y)=p^{-1}(x^p+y^p-(x+y)^p)$.
Any set theoretic map $\d:A\ra A$ from a ring $A$ to itself satisfying the above axioms will be referred to a $p$-{\it derivation}; such operators were introduced independently in  \cite{joyal,char} and they implicitly arise in the theory of Witt rings. For any such $\d$, the map $\phi:A\ra A$, $\phi(a)=a^p+p\d a$ is a ring homomorphism lifting the $p$-power Frobenius on $A/pA$; and vice versa, given a $p$-torsion free ring $A$ and a ring homomorphism $\phi:A\ra A$ lifting Frobenius the map $\d:A\ra A$, $\d a=\frac{\phi(a)-a^p}{p}$ is a $p$-derivation.

\subsubsection{Differential equations}
 Differential equations $F(t,x,\d_t x,...,\d_t^n x)=0$ (with $F$
smooth) satisfied by smooth functions $x\in C^{\infty}(\bR)$ are replaced by {\it arithmetic
differential equations} $F(\alpha,\d_p \alpha,...,\d_p^n \alpha)=0$ with $F\in R[x_0,x_1,...,x_n]\h$ satisfied by numbers
$\alpha\in R$.
As we shall see it is crucial to allow $F$ to be in the $p$-adic completion of the polynomial ring rather than in the polynomial ring itself;
indeed if one restricts to polynomial $F$'s  the main interesting examples of the theory are left out.

\subsubsection{Jet spaces} More generally, the Lie-Cartan geometric theory of differential equations has an arithmetic analogue which we explain now.  Let  $M \ra T$ be  a submersion of smooth manifolds with $\dim T=1$, and let
 $$J^n(M/T)=\{J^n_t(s);s\in \Gamma(M/T), t\in T\}$$ be the space of $n$-jets $J^n_t(s)$ of smooth sections $s$ of
$M \ra T$ at points $t \in T$. Cf. \cite{ALV, IL, olver} for references to differential geometry.
If $M=\bR^d \times \bR$, $T=\bR$,  $M\ra T$ is the second projection, and $x=(x_1,...,x_d)$,
 $t$ are global coordinates on $\bR^d$ and $\bR$  respectively then $J^n(M/T)=\bR^{(n+1)d}\times \bR$ with  {\it jet coordinates} $x,x',...,x^{(n)},t$. So
for general $M \ra T$ the map $J^n(M/T)\ra M$ is a fiber bundle with fiber $\bR^{nd}$ where $d+1=\text{dim}(M)$.
One has  the total derivative operator $$\d_t:C^{\infty}(J^n(M/T))\ra C^{\infty}(J^{n+1}(M/T))$$ which in coordinates is given by $$\d_t=\frac{\partial}{\partial t}+
\sum_{j=0}^n \sum_{i=1}^d x_i^{(j+1)}\frac{\partial}{\partial x_i^{(j)}}.$$
In the arithmetic theory the analogues of
the manifolds $J^n(M/T)$ are certain  formal schemes (called {\it arithmetic jet spaces} or {\it $p$-jet spaces}) $J^n(X)=J^n(X/R)$ defined as follows.
Assume $X$ is affine, $X=Spec\ \frac{R[x]}{(f)}$ with $x,f$ tuples; the construction that follows is easily globalized to the non-affine case. 
Let $x',...,x^{(n)},...$ be new tuples of variables, consider the polynomial ring $R\{x\}:=R[x,x',x'',...]$, 
let  $\phi:R\{x\}\ra R\{x\}$ be the unique ring homomorphism extending $\phi$ on $R$ and sending $\phi(x)=x^p+px'$,
$\phi(x')=(x')^p+px''$,..., and let   $$\d=\d_p:R\{x\}\ra R\{x\}$$ be the $p$-derivation
 $$\d F=\frac{\phi(F)-F^p}{p}.$$
Then one defines 
$$J^n(X)=Spf\ \frac{R[x,x',...,x^{(n)}]\h}{(f,\d f,...,\d^n f)}.$$

\subsubsection{Differential equations  on manifolds}
 Usual {\it differential equations} are defined geometrically  as elements of the ring $C^{\infty}(J^n(M/T))$; alternatively such elements are referred to as (time dependent) {\it Lagrangians} on $M$. Their analogue in the arithmetic theory, which we call  {\it arithmetic differential equations}  \cite{book}, are the  elements of the ring $\cO^n(X):=\cO(J^n(X))$.
 For group schemes the following concept \cite{char} plays an important role:
if $G$ is a group scheme of finite type over $R$ then we may consider the
$R$-module $\cX^n(G)=Hom_{gr}(J^n(G),\widehat{\bG}_a) \subset
\cO^n(G)$.

 Going back to arbitrary schemes of finite type $X/R$ note that
$\d_p:R\{x\}\ra R\{x\}$ induces maps $\d=\d_p:\cO^n(X)\ra \cO^{n+1}(X)$ which can be viewed as arithmetic analogues of the total derivative operator $\d_t:C^{\infty}(J^n(M/T))\ra C^{\infty}(J^{n+1}(M/T))$. The latter   is a ``generator" of the
Cartan distribution defined by $dx_i-x_i'dt$, $dx_i'-x_i''dt$, etc. Note however that the forms in differential geometry defining the Cartan distribution do not have a direct  arithmetic analogue; for one thing there is no form ``$dp$" analogous to $dt$.
On the other hand in the arithmetic case we have induced ring homomorphisms $\phi:\cO^n(X)\ra \cO^{n+1}(X)$ which have no analogue in differential geometry.

\subsubsection{Differential functions} Any differential equation $F \in C^{\infty}(J^n(M/T))$ defines a natural {\it differential function}  $F_*:\Gamma(M/T) \ra C^{\infty}(T)$; in coordinates sections $s\in \Gamma(M/T)$ correspond to functions $x=x(t)$ and then $F_*$ sends $x(t)$ into  $F_*(x(t))=F(x(t),\d_t x (t),...,\d_t^n x(t))$. Analogously any arithmetic differential equation $f \in \cO(J^n(X))$ defines a map of sets $f_*:X(R)\ra R$, referred to as a $\d$-{\it function}, which in affine coordinates sends $\alpha \in X(R)\subset R^N$ into
$f_*(\alpha):=F(\alpha,\d_p \alpha,...,\d^n_p\alpha) \in R$ if $F \in R[x,x',...,x^{(n)}]\h$ represents $f$.
If $X=G$ is in addition a group scheme and $\psi\in \cX^n(G)$  then  $\psi_*:G(R)\ra \bG_a(R)=R$ is a group homomorphism called a $\d$-{\it character} of $G$.
For $X/R$ smooth and $R/pR$  algebraically closed  $f$ is uniquely determined by $f_*$ so one can identify $f$ with $f_*$ and  $\psi$ with $\psi_*$.

\subsubsection{Prolongations of vector fields} 
For any  vertical vector field 
$$\xi:=\sum_{i=1}^d a_i(t,x)\frac{\partial}{\partial x_i}$$ 
on $M/T$ one can consider the 
canonical prolongations
 $$\xi^{(n)}:=\sum_{j=0}^n \sum_{i=1}^d (\d^j_t a_i(t,x))\frac{\partial}{\partial x_i^{(j)}}$$
 on $J^n(M/T)$. 
 The map $$\xi^{(n)}:C^{\infty}(J^n(M/T))\ra C^{\infty}(J^n(M/T))$$ is the unique $\bR$-derivation whose restriction to $C^{\infty}(M)$ is $\xi$ and which commutes with the total derivative operator $\d_t$.
 The above construction has an arithmetic analogue that plays a key technical role in the development of the theory. Indeed for any affine smooth $X/R$ and any $R$-derivation $\xi:\cO(X) \ra \cO(X)$ the canonical prolongation $$\xi^{(n)}:\cO^n(X)[1/p]\ra \cO^n(X)[1/p]$$ is defined as the unique $K=R[1/p]$-derivation whose restriction to $\cO(X)$ is $\xi$ and which commutes with $\phi$. This construction then obviously globalizes.

 \subsubsection{Infinitesimal symmetries of differential equations} Some
 Galois theoretic concepts based on  prolongations of vector fields  have  arithmetic analogues. Indeed recall that if ${\mathcal L}\subset C^{\infty}(J^n(M/T))$ is a linear subspace of differential equations then a vertical vector field $\xi$ on $M/T$ is called an {\it infinitesimal symmetry} of ${\mathcal L}$ if $\xi^{(n)}{\mathcal L}\subset {\mathcal L}$;
 it is called a {\it variational infinitesimal symmetry} if $\xi^{(n)}{\mathcal L}=0$. Similarly given an $R$-submodule  ${\mathcal L}\subset \cO^n(X)$ 
  an {\it infinitesimal symmetry} of ${\mathcal L}$ is an $R$-derivation $\xi:\cO(X)\ra \cO(X)$ such that
  $\xi^{(n)}{\mathcal L}\subset {\mathcal L}[1/p]$. One says $\xi$ is a {\it variational infinitesimal symmetry} of ${\mathcal L}$
  if $\xi^{(n)}{\mathcal L}=0$.

  \subsubsection{Total differential forms on manifolds} 
  Recall that a  {\it total differential form} (\cite{olver}, p. 351) on $M/T$
is an expression that in coordinates looks like a sum of expressions $$F(t,x,x',...,x^{(n)})dx_{j_1} \wedge...\wedge dx_{j_i},$$
  It is important to introduce an arithmetic analogue of this  which we now explain.
  We consider the case of top forms ($i=d$) which leads to what we will call {\it $\d$-line bundles}. 
Denote by $\cO^n$ the sheaf $U \mapsto \cO^n(U)$ on $X$ for the Zariski topology.
Define
a {\it $\d$-line bundle} of order $n$ on $X$  to be a locally
free $\cO^n$-module  of rank $1$. Integral powers of bundles need to be generalized as follows.
Set $W=\bZ[\phi]$ (ring of polynomials with $\bZ$-coefficients in
the symbol $\phi$). For $w=\sum a_s\phi^s$ write $deg(w)=\sum a_s$.
Also let $W_+$ be the set of all $w=\sum a_s\phi^s \in W$ with
$a_s\geq 0$ for all $s$. If $L$ is a line bundle on $X$ given by a
cocycle $(g_{ij})$ and $w=\sum_{s=0}^n a_s\phi^s\in W$, $w \neq 0$,
$a_n \neq 0$, then define a $\d$-line bundle $L^w$ of order $n$ by
the cocycle $(g_{ij}^w)$, $g_{ij}^w=\prod_s \phi^s(g_{ij})^{a_s}$.
With all these definitions in place we may define the following rings which, by the way, are the main objects of the theory:
$$R_{\d}(X,L)=\bigoplus_{0\neq w
\in W_+}H^0(X,L^w).$$
Note that the above is a graded ring without unity.
The homogeneous elements of $R_{\d}(X,L)$ can be viewed as
arithmetic analogues of Lagrangian densities and can also be
referred to as {\it arithmetic differential equations} \cite{book}.

\subsubsection{Differential forms on jet spaces and calculus of variations} The spaces $$\Lambda^i_{\uparrow}(J^n(M/T))$$  of vertical smooth $i$-forms on $J^n(M/T)$ (generated by $i$-wedge products of forms $dx,dx',...,dx^{(n)}$) play an important role in the calculus of variations \cite{olver}. These spaces fit into a deRham complex
where the differential is the vertical exterior differential $$d_{\uparrow}:\Lambda^i_{\uparrow}(J^n(M/T))\ra \Lambda^{i+1}_{\uparrow}(J^n(M/T))$$ with respect to the variables $x,x',...,x^{(n)}$. On the other hand we have  unique operators $$\d_t:\Lambda^i_{\uparrow}(J^n(M/T))\ra \Lambda^i_{\uparrow}(J^{n+1}(M/T))$$ that commute with  $d_{\uparrow}$, induce a derivation on the exterior algebra, and for $i=0$ coincide with the total derivative operators. Then one can define spaces of {\it functional forms} (\cite{olver}, p. 357)
$$\Lambda^i_*(J^n(M/T))=\frac{\Lambda_{\uparrow}^i(J^n(M/T))}{Im(\d_t)}$$
and the (vertical part of the) {\it variational complex} (\cite{olver}, p. 361) with differentials
$$d:\Lambda^i_*(J^n(M/T))\ra \Lambda^{i+1}_*(J^n(M/T)).$$
The class of $\omega \in \Lambda_{\uparrow}^i(J^n(M/T))$ in $\Lambda_*^i(J^n(M/T))$ is denoted by $\int \omega dt$. For $i=0$, $\omega=F$, have the formula $d(\int F dt)= \int EL(F)dt$ in $\Lambda_*^i(J^{2n}(M/T))$
where $EL(F)=\sum_{i=1}^d F_i dx_i$ is the {\it Euler-Lagrange} total differential  form,
$$F_i=\sum_{j=0}^n (-1)^j \d_t^j \left( \frac{\partial F}{\partial x_i^{(j)}}\right).$$
Noether's theorem then says that for any vertical vector field $\xi $ on $M/T$ we have the formula
$$\langle \xi, EL(F)\rangle - \xi^{(n)} (F)=\d_t G$$
for some $G \in C^{\infty}(J^{2n-1}(M/T))$; $G$ is unique up to a constant. If 
$\xi$ is a variational infinitesimal symmetry of $F$ then
$G$ is referred to as the {\it conservation law} attached to this symmetry; in this case if $x(t)$ is a solution of all $F_i=0$ then
$G$ evaluated at $x(t)$ will be a constant.
Analogously, for $X/R$ smooth and affine, one can consider the modules  $\Omega^i_{\cO^n(X)/R}$ (exterior powers of the inverse limit of the K\"{a}her differentials $\Omega_{\cO^n(X)\otimes (R/p^nR)/(R/p^nR)}$), and the exterior differential $$d:\Omega^i_{\cO^n(X)/R}\ra \Omega^{i+1}_{\cO^n(X)/R}.$$
Also one may consider  the operators $$\phi^*:\Omega^i_{\cO^n(X)/R}\ra \Omega^i_{\cO^{n+1}(X)/R};$$
again $\phi^*$ and $d$ commute. Also note that any element in the image of $\phi^*$ is uniquely divisible by $p^i$;
for any $\omega \in \Omega^i_{\cO^n(X)/R}$ and $r\geq 0$ we then set $\omega_r=p^{-ir}\phi^{*r} \omega$.
The operation $p^{-ir}\phi^{*r}$ is a characteristic zero version of the inverse Cartier operator.
For any element $\mu\in R$ (which we refer to as {\it eigenvalue}) one can define groups 
$$\Lambda^i_*(J^n(X))=\frac{\Omega^i_{\cO^n(X)/R}}{Im(\phi^*-\mu)}$$
that fit into a {\it variational complex} with differentials
$$d:\Lambda^i_*(J^n(X))\ra \Lambda^{i+1}_*(J^n(X)).$$
The class of $\omega \in \Omega^i_{\cO^n(X)/R}$ in $\Lambda_*^i(J^n(X))$ is denoted by $\int \omega dp$.
For $i=0$, $\mu=1$, $\omega=F\in \cO^n(X)$ we have $d (\int Fdp)=\int \{\sum_{i=1}^d F_i  \omega_n^i\}dp$ in $\Lambda^1_*(J^{2n}(X))$ where $F_i\in \cO^{2n}(X)$, $\omega^i$ is a basis of $\Omega_{\cO(X)/R}$, and  $\omega^i_n=p^{-n}\phi^{*n}\omega^i$.
Also an analogue of the Noether theorem holds in this context with  $\epsilon(F):=\sum_{i=1}^d F_i  \omega_n^i$  playing the role of the Euler-Lagrange form; indeed for any vector field $\xi$ on $X$ there exists $G \in \cO^{2n-1}(X)$ such that
$$\langle \xi^{(n)}, \epsilon(F)\rangle - \xi^{(n)} (F)=G^{\phi}-G.$$
If $\xi$ is a variational infinitesimal symmetry of $F$ then $G$ can be referred to  as a {\it conservation law}; in this case if $P\in X(R)$ is a point which is a solution to all $F_{i*}(P)=0$ then $G_*(P)^{\phi}=G_*(P)$. Part of the formalism above can be used to introduce an arithmetic analogue of the Hamiltonian formalism for which we refer to
\cite{BYM}.

\subsubsection{``Category of differential equations"}
A categorical framework for differential equations on manifolds can be introduced (cf. \cite{ALV}, for instance). 
In one variant of this the objects are locally isomorphic to
projective systems of submanifolds of $J^n(M/T)$ compatible with the total derivative operator
and morphisms are smooth maps between these, again compatible  with   the total derivative operator.
This categorical framework has an arithmetic analogue as follows. Note first that
 $\phi$ acts naturally on  $R_{\d}(X,L)$ but $\d$ does not. Nevertheless
 for any homogeneous $s \in R_{\d}(X,L)$ of degree $v$ the ring
 $R_{\d}(X,L)_{\langle s \rangle}$ of all fractions $f/s^w$ with $f$
 homogeneous of degree $wv$ has a natural $p$-derivation
 $\d$ on it.
This inductive system $X_{\d}(L)=(R_{\d}(X,L)_{\langle s \rangle},\d)_s$ of rings equipped with $p$-derivations $\d$ is an object of a natural category underlying a geometry more general than algebraic geometry which we refer to as $\d$-{\it geometry} \cite{book}.
 This geometry is an arithmetic analogue of the   categorical setting  in \cite{ALV} and also an arithmetic analogue of the Ritt-Kolchin $\d$-algebraic geometry \cite{Kolchin, hermann}.
 If one restricts to \'{e}tale maps of smooth schemes we have a functor $X \mapsto X_{\d}=X_{\d}(K^{\nu})$ from
 ``algebraic geometry"
  to our ``$\d$-geometry" by taking  $L=K^{\nu}$  to be a fixed power of the canonical  bundle; the natural choice later in the theory turns out to be the anticanonical bundle $L=K^{-1}$.
  In  $\d$-geometry $X_{\d}$ should be viewed as an infinite dimensional object.

  \subsubsection{Differential invariants} Recall the   concept   of {\it differential invariant} which plays a key role in the ``Galois theoretic" work of Lie and Cartan. Assume that  a  group $G$ acts on $M$ and $T$ such that $M \ra T$ is $G$-equivariant.
   (In this discussion assume $d:=\dim T$ is arbitrary.)   Then $G$ acts on $J^n(M/T)$ and the ring $C^{\infty}(J^n(M/T))^G$ of $G$-invariant elements in the ring $C^{\infty}(J^n(M/T))$  is called the ring of {\it differential invariants}. 
   There are two extreme cases of special interest: the case when $G$ acts trivially on $T$ and the case when  $G$ acts transitively on $T$. A special case of the first extreme case mentioned above is that in which
   $M=F \times T$, with  $G$  a Lie group acting trivially on $T$ and transitively on $F$;   this leads to the context  of Cartan's moving frame method and of Cartan connections. For $\dim T=1$, $T$ is the ``time" manifold,  $F$ is the ``physical space", and sections in $\Gamma(M/T)$ correspond to particle trajectories. A special case of the second extreme case mentioned above
   is the situation  encountered in the study of ``geometric structures" and  in the formulation of field theories, in which  $G$ is the group $\text{Diff}(T)$ of diffeomorphisms of $T$, $T$ is viewed as the ``physical space" or ``physical space-time", and $M$ is a {\it natural bundle} over $T$, i.e. a quotient 
  $M=\Gamma \backslash \text{Rep}_n(T)$
  of the bundle $\text{Rep}_n (T)\ra T$ of $n$-jets of frames $(\bR^d,0)\ra T$  by a Lie subgroup $\Gamma$ of the group 
  $\text{Aut}_n(\bR^d,0)$
  of $n$-jets  of diffeomorphisms of $(\bR^d,0)$; cf. \cite{ALV}, pp. 150-153 and 183. 
  (E.g. the Riemannian metrics on $T$ identify with the sections in $\Gamma(M_{1,O(d)}/T)$ where $O(d) < GL(d)=\text{Aut}_1(\bR^d,0)$ is the orthogonal group.)
  In the  case when the action of $G$ on $T$ is trivial the ring of differential invariants above  turns out to have an interesting arithmetic analogue, namely the  ring of {\it $\d$-invariants of a correspondence}; cf. the discussion below.
  The special case of Cartan connections has also an arithmetic analogue: the {\it arithmetic logarithmic derivative} attached to a $\d$-{\it flow}. The latter has   a flavor different from that  of $\d$-invariants of correspondences and  will be discussed later. 
 There are also interesting candidates for arithmetic analogues  of the situation when $G=\text{Diff}(T)$; cf. our discussion of Lie groupoids below.

We next explain the ring of $\d$-invariants of a correspondence.
Pursuing an analogy with usual geometric invariant theory assume we are
given $(X,L)$,  a {\it correspondence} on $X$ (i.e. a morphism $\sigma=(\sigma_1,\sigma_2):\tilde{X} \ra X \times X$),
and a {\it linearization} of $L$ (i.e. an isomorphism $\beta:\sigma_1^*
L \simeq \sigma_2^*L$). Then we may  define the ring of {\it $\d$-invariants of
$\sigma$} by $$R_{\d}(X,L)^{\sigma}=\{f\in R_{\d}(X,L);\ \beta
\sigma_1^*f=\sigma_2^*f\}.$$ (This ring, again, has no unity!)
The homogeneous elements $s$ of $R_{\d}(X,L)^{\sigma}$ can be viewed as
arithmetic analogues of total differential forms  invariant under
appropriate symmetries. The  inductive system of rings $$(R_{\d}(X,L)^{\sigma}_{\langle s \rangle},\d)$$ equipped with $p$-derivations  $\d$  can be viewed as an incarnation of the quotient space
 ``$X_{\d}/\sigma_{\d}$" in $\d$-geometry (and, under quite general hypotheses, is
 indeed the categorical quotient of $X_{\d}$ by $\sigma_{\d}$); note that, in most
 interesting examples (like the ones in Theorem \ref{thm2} below),
  the  quotient $X/\sigma$ does not
 exist in usual algebraic geometry (or rather the categorical
  quotient in algebraic geometry reduces to a point).
  
 It is worth revisiting the sequence of ideas around differential invariants. Classical Galois theory deals with  algebraic equations (satisfied by numbers).
Lie and Cartan extended Galois' ideas, especially through the concept of differential invariant, to the study of differential equations (satisfied by functions); roughly speaking they replaced numbers by functions. In the theory  presented here functions are replaced back by numbers. But we did {\it not} come back to where things started because we have added new structure, the operator $\d=\d_p$.  In particular the $\d$-invariants mentioned above, although arithmetic in nature, and although attached to algebraic equations (defining $X$, $\tilde{X}$, $\sigma$), are nevertheless {\it not} ``Galois theoretic" in any classical sense.

 \subsubsection{Lie groupoids}  Let us go back to the (interrelated) problems
   of finding arithmetic analogues of $T\times T$ and of $\text{Diff}(T)$.
    The arithmetic analogue of  $S:=T\times T$    is usually referred to as the hypothetical tensor product ``$\bZ \otimes_{\bF_1}\bZ$" over the ``field with one element"; cf. \cite{maninf1} for some history of this. 
    The arithmetic analogue of $\text{Diff}(T)$ could be thought of as the Galois group ``$\text{Gal}(\bQ/\bF_1)$".
    According to a suggestion of Borger \cite{borger} one should take ``$\bZ \otimes_{\bF_1}\bZ$" to be, by definition, the big Witt ring of $\bZ$. Since we are here in a local arithmetic situation it is convenient, for our purposes, to take, as an arithmetic analogue of $S=T \times T$, the schemes $\Sigma_m=Spec\ W_m(R)$, where $W_m$ is the $p$-typical functor
   of Witt vectors of length $m+1$ (we use Borger's indexing) and $R$ is a complete local $p$-ring.   
  On the other hand an infinitesimal analogue of $\text{Diff}(T)$ is the Lie groupoid of $T$ defined as the projective system
  of groupoids ${\mathcal G}_n(T):=J^n(S/T)^*$ (where the upper $*$ means taking ``invertible" elements and $S=T \times T \ra T$ is the second projection.) So an arithmetic analogue of the system $J^n(S/T)^*$ (which at the same time would be an infinitesimal analogue of ``$\text{Gal}(\bQ/\bF_1)$") would be the system $J^n(\Sigma_m)$ equipped with the natural maps induced by the comonad map. This system is a rather non-trivial object \cite{pw}.  It is worth noting that  a good arithmetic analogue of $J^n(S/T)^*$ 
    could also be the ``usual" jet spaces (in the sense of \cite{annals}) of $W_m(R)/R$ (which in this case
   can be constructed directly from the module of
     K\"{a}her differentials $\Omega_{W_m(R)/R}$). By the way  $\Omega_{W_m(R)/R}$ is  also  the starting point for  the construction of the deRham-Wit complex \cite{he}. However $\Omega$ involves usual derivations (rather than Fermat quotients)  so taking $\Omega$ as a path to an arithmetic analogue of $J^n(S/T)^*$  seems, again, like ``going arithmetic halfway". 
     On the contrary taking the system $J^n(\Sigma_m)$ as the analogue of $J^n(S/T)^*$  seems to achieve, in some sense, the task of
     ``going arithmetic all the way".
     
     We end by remarking that if we denote $\pi_1$ and $\pi_2$ the source and target  projections from ${\mathcal G}_n(T)$ into $T$ then
     there are natural ``actions" $\rho$ of the groupoids ${\mathcal G}_n(T)$ on all natural bundles 
     $M_{n,\Gamma}:=\Gamma \backslash \text{Rep}_n(T) \ra T$
     fitting into diagrams:
     $$
     \begin{array}{rcl}
     M_{n,\Gamma} \times_{T,\pi_1} {\mathcal G}_n(T)  & \stackrel{\rho}{\longrightarrow} & M_{n,\Gamma} \\
     p_2 \downarrow & \  & \downarrow\\
     {\mathcal G}_n(T) & \stackrel{\pi_2}{\longrightarrow} & T
     \end{array}
     $$
     where $p_2$ is the second projection. The above induce ``actions" 
      $$
     \begin{array}{rcl}
     J^n(M_{n,\Gamma}/T) \times_{T,\pi_1} {\mathcal G}_n(T)  & \stackrel{\rho}{\longrightarrow} & J^n(M_{n,\Gamma}/T) \\
     p_2 \downarrow & \  & \downarrow \\
     {\mathcal G}_n(T) & \stackrel{\pi_2}{\longrightarrow} & T
     \end{array}
     $$
     where $p_2$ is again the second projection.
One can consider rings of differential invariants 
$$C^{\infty}(J^n(M_{n,\Gamma}/T))^{\rho}:=\{F\in C^{\infty}(J^n(M_{n,\Gamma}/T));F \circ \rho=F \circ p_1\}$$
where $p_1:J^n(M_{n,\Gamma}/T) \times_{T,\pi_1} {\mathcal G}_n(T)\ra J^n(M_{n,\Gamma}/T)$ is the fist projection.
     There should be arithmetic analogues of the above ``actions" and rings of differential invariants. One can argue that the analogue of $\text{Rep}_n(T)$ is, again,
     the system $J^n(\Sigma_m)$. Then for $\Gamma=1$ the analogue of $J^n(M_{n,\Gamma}/T)$ might be $J^n(J^n(\Sigma_m))$; the analogue, for $\Gamma=1$,
     of the ``action" $\rho$ above could then be the map $J^n(J^n(\Sigma_{m',m''}))\ra J^n(J^n(\Sigma_{m'+m''}))$ where $\Sigma_{m',m''}=Spec\ W_{m'}(W_{m''}(R))$.
     A challenge would then be to find arithmetic analogues of non-trivial $\Gamma$s.

\subsubsection{Differential Galois theory of linear equations}
This subject is best explained in a complex (rather than real) situation.
Classically (following Picard-Vessiot and Kolchin \cite{kolchin}) one starts with a differential field ${\mathcal F}$ of meromorphic functions on a domain $D$ in the complex plane ${\mathbb C}$ and one fixes an $n \times n$ matrix
$A\in {\mathfrak g}{\mathfrak l}_n({\mathcal F})$ in the Lie algebra of the general linear group
$GL_n({\mathcal F})$. The problem then is to develop a ``differential Galois theory" for  equations of the form
$$\d_z U=A\cdot U$$
where $U\in GL_n({\mathcal G})$, ${\mathcal G}$ a field of meromorphic functions on a subdomain of $D$, and $z$ a coordinate on $D$. The start of the theory is as follows. One fixes a solution $U$ and introduces the differential Galois group $G_{U/{\mathcal F}}$ of $U/{\mathcal F}$ as the group of all 
${\mathcal F}-$automorphisms of the field ${\mathcal F}(U)$ that commute with $d/dz$. One can ask for an arithmetic analogue of this theory. There is a well developed difference algebra analogue of linear differential equations \cite{SVdP}; but the arithmetic differential theory is still in its infancy \cite{adel1,adel2,adel3}. What is being proposed in loc. cit. in the arithmetic theory is to fix a matrix $\alpha\in {\mathfrak g}{\mathfrak l}_n(R)$ and define $\d$-linear equations as equations of the form
$$\d u=\alpha \cdot u^{(p)}$$
where $u=(u_{ij})\in GL_n(R)$, $\d u:=(\d u_{ij})$,  and $u^{(p)}:=(u_{ij}^p)$. 
Note that the above equation 
is equivalent to $\phi(u)=\epsilon \cdot u^{(p)}$ where $\epsilon=1+p\alpha$ which 
is  {\it not} a  difference equation for $\phi$ in the sense of \cite{SVdP}; indeed difference equations for $\phi$ have the form $\phi(u)=\epsilon \cdot u$. 
To define the $\d$-Galois group of such an equation start with a  $\d$-subring $\cO\subset R$ and let $u\in GL_n(R)$
be a solution of our equation. Let $\cO[u]\subset R$ the ring generated by the entries of $u$; clearly $\cO[u]$ is a $\d$-subring of $R$. We define the {\it $\d$-Galois group} $G_{u/\cO}$ of $u/\cO$
as the subgroup of all $c\in GL_n(\cO)$ for which there exists an $\cO$-algebra automorphism
$\sigma$ of $\cO[u]$ such that $\sigma\circ \d=\d\circ \sigma$ on $\cO[u]$ and such that $\sigma(u)=uc$. The theory starts from here.

\subsubsection{Cartan connections} The Maurer-Cartan connection attached to a Lie group $G$ is a canonical map $T(G)\ra L(G)$ from the tangent bundle $T(G)$ to the Lie algebra $L(G)$; for $G=GL_n$ it is given by the form
$dg\cdot g^{-1}$ and its algebraic incarnation is Kolchin's logarithmic derivative map \cite{kolchin}
$GL_n({\mathcal G})\ra {\mathfrak g}{\mathfrak l}_n({\mathcal G})$, $u\mapsto \d_z u\cdot u^{-1}$. In our discussion of linear equations above the arithmetic analogue of the Kolchin logarithmic derivative map is
 the map $GL_n(R)\ra {\mathfrak g}{\mathfrak l}_n(R)$,
$$u\mapsto \d u \cdot (u^{(p)})^{-1}.$$
 This map is naturally attached to the lift of Frobenius
$\phi_{GL_n,0}:\widehat{GL_n}\ra \widehat{GL_n}$ whose effect on the ring
$\cO(\widehat{GL_n})=R[x,\det(x)^{-1}]\h$ is $\phi_{GL_n,0}(x)=x^{(p)}$.
The latter lift of Frobenius behaves well (in a precise sense to be explained later) with respect to the maximal torus $T\subset GL_n$ of diagonal matrices and with respect to the Weyl group $W\subset GL_n$ of permutation matrices but it behaves ``badly" with respect to other subgroups like the classical groups $SL_n,SO_n,Sp_n$.
(This phenomenon does not occur in the geometry of Lie groups where the Maurer-Cartan form behaves well with respect to {\it any} Lie subgroup of $GL_n$, in particular with respect to the classical groups.)
In order to remedy the situation one is lead to generalize the above constructions by replacing $\phi_{GL_n,0}$ with other lifts of Frobenius
$\phi_{GL_n}:\widehat{GL_n}\ra \widehat{GL_n}$ that are adapted to each of these classical groups. Here is a description of the resulting framework.

First we define an arithmetic analogue of the Lie algebra ${\mathfrak g}{\mathfrak l}_n$ of $GL_n$ as the set ${\mathfrak g}{\mathfrak l}_n$  of $n\times n$ matrices over $R$ equipped with the non-commutative group law $+_{\d}:{\mathfrak g}{\mathfrak l}_n\times {\mathfrak g}{\mathfrak l}_n\ra {\mathfrak g}{\mathfrak l}_n$ given by
$$a+_{\d}b:=a+b+pab,$$
where the addition and multiplication in the right hand side are those of ${\mathfrak g}{\mathfrak l}_n$, viewed as an associative algebra.
There is a natural ``$\d$-adjoint" action $\star_{\d}$ of $GL_n$ on ${\mathfrak g}{\mathfrak l}_n$ given by
$$a\star_{\d} b:=\phi(a) \cdot b \cdot \phi(a)^{-1}.$$

Assume now one is given a ring endomorphism $\phi_{GL_n}$ of $\cO(\widehat{GL_n})$ lifting Frobenius, i.e. a ring endomorphism whose reduction mod $p$ is the $p$-power Frobenius on 
$\cO(GL_n)\h/(p)=k[x,\det(x)^{-1}]$; we still denote by $\phi_{GL_n}:\widehat{GL_n}\ra \widehat{GL_n}$ the induced morphism of $p$-formal schemes and we refer to it as a lift of Frobenius on $\widehat{GL_n}$ or simply as a $\d$-{\it flow} on $\widehat{GL_n}$.
Consider the matrices $\Phi(x)=(\phi_{GL_n}(x_{ij}))$ and $x^{(p)}=(x_{ij}^p)$
with entries in $\cO(GL_n)\h$; then
 $\Phi(x)=x^{(p)}+p\Delta(x)$ where $\Delta(x)$ is some matrix with entries in $\cO(GL_n)\h$. Furthermore given a lift of Frobenius $\phi_{GL_n}$ as above one defines, as usual, a $p$-derivation
$\d_{GL_n}$ on $\cO(GL_n)\h$ by setting $\d_{GL_n}(f)=\frac{\phi_{GL_n}(f)-f^p}{p}.$

 Assume furthermore that we are given a smooth closed subgroup scheme $G\subset GL_n$. We say that $G$ is $\phi_{GL_n}$-horizontal if $\phi_{GL_n}$ sends the ideal of $G$ into itself; in this case we have a lift of Frobenius endomorphism $\phi_G$ on $\widehat{G}$, equivalently on $\cO(G)\h$. 

Assume the ideal of $G$ in $\cO(GL_n)$ is generated by polynomials $f_i(x)$.
Then recall that the Lie algebra $L(G)$ of $G$  identifies, as an additive group, to the  subgroup of the usual additive group $({\mathfrak g}{\mathfrak l}_n,+)$ consisting of all matrices $a$ satisfying
 $$``\epsilon^{-1}"f_i(1+\epsilon a)=0,$$
where $\epsilon^2=0$. Let
 $f^{(\phi)}_i$ be the polynomials obtained from $f_i$ by applying  $\phi$ to the coefficients.
Then 
we  define the $\d$-Lie algebra $L_{\d}(G)$ as the  subgroup of $({\mathfrak g}{\mathfrak l}_n,+_{\d})$ consisting of all the matrices $a\in {\mathfrak g}{\mathfrak l}_n$ satisfying
$$p^{-1}f_i^{(\phi)}(1+pa)=0.$$
 
 The analogue of Kolchin's logarithmic derivative (or of the Maurer-Cartan connection)  will then be
the map, referred to as the {\it arithmetic logarithmic derivative}, $l\d:GL_n\ra {\mathfrak g}{\mathfrak l}_n,$
defined by
$$l\d a:=\frac{1}{p}\left(\phi(a)\Phi(a)^{-1}-1\right)=(\d a-\Delta(a))(a^{(p)}+p\Delta(a))^{-1}.$$
For $G$ a $\phi_{GL_n}$-horizontal subgroup $l\d$ above induces a map 
$l\d:G\ra L_{\d}(G)$.
Now  any $\alpha\in L_{\d}(G)$  defines  an equation of the form
$
l\d u=\alpha,
$
with unknown $u\in G$; such an equation will be referred to as a $\d$-linear equation (with respect to our $\d$-flow). Later in the paper we will explain our results about existence of $\d$-flows on $GL_n$ compatible with the classical groups. These $\d$-flows will produce, as explained above, corresponding $\d$-linear equations.

\subsection{Main task of the theory}
At this point we may formulate  the main technical tasks of the theory.
Let, from
now on, $R=\widehat{\bZ_p^{ur}}$.
First
given a specific scheme $X$ (or group scheme $G$) the task is to compute the
rings $\cO^n(X)$ (respectively the modules $\cX^n(G)$).
More generally given a specific pair $(X,L)$ we want  to compute the ring $R_{\d}(X,L)$.
Finally given $(X,L)$,  a correspondence $\sigma$ on $X$, and a
linearization of $L$,
we want to  compute the ring $R_{\d}(X,L)^{\sigma}$.
 {\it The main applications of the theory (cf. the subsection below on motivations)
 arise as a result of the presence of interesting/unexpected elements and relations in the rings $\cO^n(X)$, $R_{\d}(X,L)$, $R_{\d}(X,L)^{\sigma}$.}

\subsection{Motivations of the theory}
There are a number of motivations for developing such a
theory.

\subsubsection{ Diophantine geometry}
   Usual differential equations satisfied by functions can be used to
prove diophantine results over function fields (e.g. the function
field analogues of the Mordell conjecture \cite{manin} and of the
Lang conjecture \cite{annals}). In the same vein one can hope to use
``arithmetic differential equations'' satisfied by numbers to prove
diophantine results over number fields. This idea actually works in
certain situations, as we will explain below. Cf. \cite{pjets,local}.
The general strategy is as follows. Assume one wants to prove that a set
$S$ of points on an algebraic variety is finite. What one tries to do is
find a system of arithmetic differential equations $F_i(\alpha,\d_p \alpha,...,\d_p^n \alpha)=0$
satisfied by all  $\alpha \in S$; then one tries,
 using algebraic operations and the ``differentiation" $\d_p$,
 to eliminate  $\d_p \alpha,...,\d_p^n \alpha$
from this system to  obtain another system $G_j(\alpha)=0$ satisfied by all $\alpha \in S$, where $G_j$ do not involve the ``derivatives" anymore.
Finally one proves, using usual algebraic geometry that the latter system has only finitely many solutions.

\subsubsection{ ``Impossible" quotient spaces} If $X$ is an algebraic
variety and $\sigma:\tilde{X} \ra X\times X$ is a correspondence on $X$ then the categorical quotient
$X/\sigma$ usually reduces to a point in the category of algebraic varieties. In some
sense this is an illustration of the limitations of classical
algebraic geometry and suggests the challenge of creating more
general geometries in which $X/\sigma$ becomes interesting. (The
non-commutative geometry of A. Connes \cite{Connes} serves in particular this
purpose.)
As explained above we proposed, in our work, to replace the algebraic equations
of usual algebraic geometry by ``arithmetic differential equations";
the resulting new geometry is called $\d_p$-geometry (or simply
$\d$-geometry). Then it turns out that important class of quotients
$X/\sigma$ that reduce to a point in usual algebraic geometry become
interesting in $\d$-geometry (due to the existence of interesting $\d$-invariants).
 A general principle seems to emerge
according to which this class coincides with the class of
``analytically uniformizable" correspondences. Cf. \cite{book}.

\subsubsection{ ``Impossible" liftings to characteristic zero} A series of phenomena
belonging to algebraic geometry in characteristic $p$, which do not
lift to characteristic $0$ in algebraic geometry, can be lifted
nevertheless to  characteristic $0$ in $\d$-geometry.
This seems to be a quite general principle with various incarnations
throughout the theory (cf. \cite{difmod,book,igusa}) and
illustrates, again, how a limitation of classical algebraic geometry
can be overcome by passing to $\d$-geometry.

 \subsection{Comparison with other theories} It is worth noting that
the paradigm of our work
 is quite different  from the following other paradigms:

 1) Dwork's  theory of
  $p$-adic differential equations
\cite{dwork} (which is a theory about $\delta_t$ acting on functions
in $\bQ_p[[t]]$ and not about $\d_p$ acting on numbers; also Dwork's theory is a theory of linear differential equations
whereas the theory here is about analogues of non-linear differential equations),

 2) Vojta's
jet spaces \cite{vojta} (which, again, although  designed for
arithmetic purposes, are nevertheless constructed from Hasse-Schmidt
derivations ``$\frac{1}{n!}\delta_t^n$" acting on functions and not
from operators acting on numbers),

 3) Ihara's differentiation of
integers \cite{I} (which, although based on Fermat quotients,  goes from characteristic zero to
characteristic $p$ and hence, unlike our $\d_p$, cannot be iterated),

 4)
the point of view of Kurokawa et. al. \cite{kur} (which uses an
operator on numbers very different from $\d_p$ namely
$\frac{\partial \alpha}{\partial p}:=np^{n-1}\beta$ for $\alpha=p^n
\beta \in {\mathbb Z}$, $p\not| \beta$),

5) the theory of Drinfeld modules \cite{Drinfeld}(which is entirely
in characteristic $p$),

6) the difference geometry in the work of Cohn,
Hrushovski-Chatzidakis \cite{CH}, and others (in which the jet
spaces are $n$-fold products of the original varieties as opposed to
the jet spaces in our work which are, as we shall see,  bundles over the original
varieties with fibers affine spaces),

7) Raynaud's deformation to Witt vectors $W_2(k)$ over a field $k$
of characteristic $p$ \cite{raynaud} (which again leads to operators
from characteristic zero to characteristic $p$ which cannot be
iterated; loosely speaking $W_2(k)$ in Raynaud's approach is
replaced in our theory by $W_2(W(k))$).

8) the work of Soul\'{e}, Deitmar, Connes, and many others on the
``geometry over the field $\bF_1$ with one element". (In their work
passing from the geometry over $\bZ$ to the geometry over $\bF_1$ amounts to {\it removing} part of the structure defining commutative rings, e.g. removing
 addition and hence considering multiplicative monoids instead of rings. On the contrary  our approach (cf. the Introduction to
\cite{book}), and  Borger's [3], can be seen as a way of passing from $\bZ$ to $\bF_1$   by {\it adding} structure to the commutative rings, specifically adding the operator(s) $\d_p$).

9) the work of Joyal \cite{joyal}  on the Witt functor (which is a right adjoint to the forgetful functor from ``$\d$-rings" to rings as opposed to our arithmetic jet  functor which is a left adjoint to the same forgetful functor. As it is usually the case the left and right stories turn out to be rather different.)

10) the theory of the Greenberg transform, cf. Lang's thesis  and  \cite{greenberg} (which attaches to a scheme $X/R$  varieties $G^n(X)$ over $k$; one can show \cite{pjets}  that $G^n(X)\simeq J^n(X)\otimes_R k$ so the arithmetic jet spaces  are certain remarkable liftings to characteristic zero of the Greenberg transforms. The operators $\d$ on $\cO^n(X)$ do not survive after reduction mod $p$ as operators on the Greenberg transforms.)

11) the work on the deRham-Witt complex, cf., e.g. \cite{he} (which has as its starting point the study of K\"{a}hler differential of Witt vectors; on the contrary, what our research suggests, cf. \cite{pw}, is to push arithmetization one step further by analyzing instead the {\it arithmetic} jet spaces of Witt vectors.)

12) the theory $\phi$-modules (which is a theory about linear equations as opposed to the theory here which is non-linear).

\section{Main results}
We present in what follows a sample of our  results. We always set   $\overline{A}=A \otimes_{\mathbb Z}\bF_p=A/pA$, $\overline{X}=X \otimes_{\mathbb Z} \bF_p$ for any ring $A$ and any scheme $X$ respectively. Recall that we denote by $\widehat{A}$ the $p$-adic completion of $A$; for $X$ Noetherian we denote by $\widehat{X}$ the $p$-adic completion of $X$. Also, in what follows,  $R:=\widehat{\bZ_p^{ur}}$, $k:=R/pR$.
 We begin with completely general facts:

\subsection{Affine fibration structure of $p$-jet spaces}

 \begin{theorem}
 \label{fibb} \cite{char} \ 
 
1) If $X/R$ is a smooth scheme of relative dimension $d$ then $X$ has an affine covering $X_i$ such that
$J^n(X_i)\simeq \widehat{X_i} \times \widehat{\bA^{nd}}$ in the category of $p$-adic formal schemes.

2) If $G/R$ is a smooth group scheme of relative dimension $d$,
with formal group law ${\mathcal F}\in R[[T_1,T_2]]^d$ ($T_1,T_2$  $d$-tuples), then the kernel
of $J^n(G)\ra \widehat{G}$ is isomorphic to $\widehat{\bA^{nd}}$ with composition law obtained from the formal series 
$$\d {\mathcal F},...,\d^n {\mathcal F}\in R[[T_1,T_2,...,T_1^{(n)},T_2^{(n)}]]^d$$ by setting $T_1=T_2=0$. 
\end{theorem}

Note that after setting $T_1=T_2=0$ the series $\d^n {\mathcal F}$ become restricted i.e. elements of $(R[T_1',T_2',...,T_1^{(n)},T_2^{(n)}]\h)^d$
so they define morphisms in the category of $p$-adic formal schemes. Assertion 1) in the  theorem makes $p$-jet spaces resemble the usual jet spaces of the Lie-Cartan theory. Note however that, even if $G$ is commutative, the kernel
of $J^n(G)\ra \widehat{G}$ is not, in general the additive group raised to some power.
Here is the idea of the proof of 1) for $n=1$. We may assume $X=Spec\ A$, $A=R[x]/(f)$,  and there is an \'{e}tale map $R[T]\subset A$ with $T$ a $d$-tuple of indeterminates. Consider the unique ring homomorphism $R[T]\ra W_1(A[T'])$ sending $T\mapsto (T,T')$ where $W_1$ is the functor of Witt vectors of length $2$ and $T'$ is a $d$-tuple of indeterminates. Using the fact that $R[T]\subset A$ is formally \'{e}tale and the fact that
the first projection  $W_1(A[T']/(p^n))\ra A[T']/(p^n)$ has a nilpotent kernel,
one constructs  a compatible sequence of homomorphisms $A \ra W_1(A[T']/(p^n))$, $T\mapsto (T,T')$. Hence one gets
a homomorphism $A\ra W_1(A[T']\h)$, 
$a\mapsto (a,\d a)$. Then one defines a homomorphism
 $R[x,x']\h/(f,\d f) \ra A[T']\h$ by sending $x \mapsto a:=\text{class}(x)\in A$, $x'\mapsto \d a$. Conversely one defines a homomorphism $A[T']\h\ra R[x,x']\h/(f,\d f)$ by sending $T'\mapsto \d T$. The two homomorphisms are inverse to each other
 which ends the proof of 1) for $n=1$.

\subsection{$\d$-functions and $\d$-characters on  curves}
The behavior of the rings $\cO^n(X)$ for smooth projective curves $X$ depends on the genus of $X$ as follows:

\begin{theorem}
\label{thm1} \cite{char,pjets,je} Let $X$ be a smooth projective curve
over $R$ of genus $g$.

1) If $g=0$ then $\cO^n(X)=R$ for all $n \geq 0$.

2) Let $g=1$. If $X$ is not a canonical lift then $\cO^1(X)=R$ (hence $\cX^1(X)=0$)  and
$\cX^2(X)$ is a free module of rank $1$; in particular $\cO^2(X)\neq
R$. 
If, on the other hand, $X$ is a canonical lift then $\cO^1(X)=\cO(\widehat{{\mathbb A}^1})$ and $\cX^1(X)$ is free of rank one.

3) If $g \geq 2$ then $J^n(X)$ is affine for $n \geq 1$; in
particular $\cO^1(X)$ separates the points of $X(R)$.
\end{theorem}

The proof of 1) is a direct computation. The idea of proof of the statements about $\cX^n$ in 2) is as follows.  Let
$N^n=Ker(J^n(X)\ra \widehat{X})$. Then one first proves (using Theorem \ref{fibb})
that $Hom(N^n,\widehat{\mathbb G}_a)$ 
has rank at least $n$ over $R$ and one computes the ranks of $\cX^n(X)$ by looking at the exact squence
$$Hom(J^2(X),\widehat{\mathbb G}_a)\ra Hom(N^2,\widehat{\mathbb G}_a)\ra H^1(X,\cO).$$
Here $Hom=Hom_{gr}$. 
The proof of 3) is based on representing $\overline{J^1(X)}$ as 
${\mathbb P}({\mathcal E})\backslash D$ where 
 ${\mathcal E}$ is a rank $2$ vector bundle on $\overline{X}$, and $D$ is an ample divisor.

If $g=1$ and $X$ is not a canonical lift then a basis $\psi$ for
$\cX^2(X)$ can be viewed as an analogue of the ``Manin map" of an
elliptic fibration \cite{manin}. Also note that
assertion 3) in Theorem \ref{thm1} implies the effective version of
the Manin-Mumford conjecture \cite{pjets}. Indeed Manin and Mumford
conjectured that if $X$ is a complex curve of genus $\geq 2$
embedded into its Jacobian $A$ then $X \cap A_{tors}$ is a finite
set. This was proved by Raynaud \cite{raynaud}. Mazur later asked
\cite{Mazur} if $\sharp(X \cap A_{tors})\leq C(g)$ where $C(g)$ is a
constant that depends only on the genus $g$ of $X$. Using 3) in
Theorem \ref{thm1} one can prove:

\begin{theorem}
\label{zuzu} \cite{pjets} For a smooth projective complex curve $X$
in its Jacobian $A$ we have
 $\sharp(X \cap A_{tors})\leq C(g,p)$ where $C(g,p)$ is
a constant that depends only on the genus $g$  and (in case $X$ is
defined over $\overline{\bQ}$) on the smallest prime $p$ of good
reduction of $X$.
\end{theorem}

Roughly speaking the idea of proof is as follows. First one can replace the complex numbers by $R$
and $A_{tors}$ by its prime to $p$ torsion subgroup $\Gamma < A(R)$. Then one embeds $X(R)\cap \Gamma$ (via the ``jet map")
into the the set of $k$-points of $\overline{J^1(X)}\cap p \overline{J^1(A)}$. But the latter is a finite set because $\overline{J^1(X)}$ is affine, $p \overline{J^1(A)}$ is projective, and both are closed in $\overline{J^1(A)}$. Moreover the  cardinality of this finite set can be bounded using B\'{e}zout's theorem.

One can ask for  global vector fields on a smooth projective curve  $X/R$ that are infinitesimal symmetries for given $\d$-functions on $X$. The only non-trivial care is that of genus $1$ (elliptic curves); indeed for genus $0$ there are no non-constant $\d$-functions (cf. Theorem \ref{thm1}) while for genus $\geq 2$ there are no non-zero vector fields.
Here is the result:

\begin{theorem}
\cite{char,book}
Let $X$ be an elliptic curve over $R$ with ordinary reduction. 

1) If $X$ has  Serre-Tate parameter $q(X)\not\equiv 1$ mod $p^2$ then there exists a non-zero global vector field on $X$ which is a variational infinitesimal symmetry for all the modules $\cX^n(X)$, $n \geq 1$.

2) If $X$ is a canonical lift (equivalently has Serre-Tate parameter $q(X)=1$) then there is no non-zero global vector field on $X$ which is a variational  infinitesimal symmetry of $\cX^1(X)$.
\end{theorem}

Finally here is a computation of differentials of $\d$-characters.

\begin{theorem}
\cite{from,book}
Assume $X$ is an elliptic curve over $R$ which is not a canonical lift and comes from an elliptic curve $X_{\bZ_p}$
over $\bZ_p$. Let $a_p\in \bZ$ be the trace of Frobenius of the reduction mod $p$ of $X_{\bZ_p}$ and let $\omega$ be a basis for the global $1$-forms on $X_{\bZ_p}$. Then there exists an $R$-basis $\psi$ of $\cX^2(X)$ such that
$$p\cdot d \psi=(\phi^{*2}-a_p\phi^*+p)\omega.$$
In particular for the eigenvalue $\mu=1$ we have
$$p\cdot d\left( \int \psi dp\right)=(1-a_p+p)\cdot \int \omega dp.$$
\end{theorem}

A similar result holds in case $X$ is a canonical lift.

\subsection{$\d$-invariants of correspondences}
Here is now a (rather roughly stated) result about $\d$-invariants
of correspondences on curves; for precise statements we refer to \cite{book}.

\begin{theorem}
\label{thm2} \cite{book} The ring
$R_{\d}(X,K^{-1})^{\sigma}$ is ``$\d$-birationally equivalent" to the ring 
$R_{\d}(\bP^1,\cO(1))$ if the correspondence $\sigma$ on $X$ ``comes
from" one of the following cases:

1) (spherical case) The standard action of $SL_2(\bZ_p)$ on $\bP^1$.

2) (flat case) A dynamical system $\bP^1 \ra \bP^1$ which is
post-critically finite with (orbifold) Euler characteristic zero.

3) (hyperbolic case) The action of  a Hecke correspondence on a
modular (or Shimura) curve.
\end{theorem}

Here by saying that $\sigma$ ``comes from" a group action on $X$ (respectively from an endomorphism of $X$) we mean that
(``up to some specific finite subschemes" for which we refer to \cite{book})
 $\tilde{X}$ is the disjoint
union of  the graphs of finitely many automorphisms generating the action
(respectively $\tilde{X}$ is the graph of the endomorphism).
Also $\d$-birational equivalence means isomorphism (compatible with
the actions of $\d$) between the $p$-adic completions of the rings
of homogeneous fractions of degree zero with denominators not
divisible by $p$. 
The proofs behind the spherical case involve direct computations. The proofs behind the flat case
use the arithmetic Manin map, i.e. the space ${\mathcal X}^2(X)$ in Theorem \ref{thm1}, plus an induction in which canonical prolongations of vector fields play a crucial role.
The proofs behind the hyperbolic case of Theorem \ref{thm2} are based on a theory of
{\it $\d$-modular forms} which we quickly survey next.

\subsection{$\d$-modular forms}
Cf. \cite{difmod,Barcau,book,hecke,eigen,igusa}. 
Let $X_1(N)$ be the complete modular curve over $R$ of level $N>4$
and let $L_1(N)\ra X_1(N)$ be the line bundle with the property that
the sections of its various powers are the classical modular forms
on $\Gamma_1(N)$ of various weights. Let $X$ be $X_1(N)$ minus the
cusps and the supersingular locus (zero locus of the Eisenstein form $E_{p-1}$ of weight $p-1$). Let
$L\ra X$ be the restriction of the above line bundle and let $V$ be
$L$ with the zero section removed. So $L^2\simeq K$. The elements of
$M^n=\cO^n(V)$ are called $\d$-modular functions. Set
$M^{\infty}=\cup M^n$. The elements of $\cO^n(X)$ are called
$\d$-modular forms of weight $0$. For $w\in W$, $w\neq 0$, the
elements of $H^0(X,L^w)\subset M^{\infty}$ are called $\d$-modular
forms of weight $w$. We let $\sigma=(\sigma_1,\sigma_2):\tilde{X} \ra X \times X$ be the union of all
 the (prime to $p$) Hecke correspondences.
 Any $\d$-modular function  $f\in M^n$ has
a ``$\d$-Fourier expansion" in $R((q))[q',...,q^{(n)}]\h$; setting
$q'=q''=...=0$ in this $\d$-Fourier expansion one gets a series in
$R((q))\h$ called the Fourier expansion of $f$.
Finally let $a_4$ and $a_6$ be variables; then consideration of the elliptic curve
$y^2=x^3+a_4 x+a_6$ yields an $R$-algebra map 
$$R[a_4,a_6,\Delta^{-1}]\ra M^0$$
where $\Delta=\Delta(a_4,a_6)$ is the discriminant polynomial. By universality we have  induced
homomorphisms
$$R[a_4,a_6,...,a_4^{(n)},a_6^{(n)},\Delta^{-1}]\h \ra M^n$$
that are compatible with $\d$.

\begin{example}
 \label{cuc}
 \cite{eigen,local}
 \textup{For $f=\sum a_nq^n$ a classical newform
over $\bZ$ of weight $2$, we get a $\d$-modular form of weight $0$
and order $2$
$$f^{\sharp}:J^2(X) \subset J^2(X_1(N)) \stackrel{J^2(\Phi)}{\ra} J^2(E_f) \stackrel{\psi}{\ra}
\widehat{\bG_a}$$ where $\Phi:X_1(N)\ra E_f$ is the Eichler-Shimura map to the
corresponding elliptic curve $E_f$ (assumed for simplicity to be
non-CM), and $\psi$ the ``unique" $\d$-character of order $2$. Then
the $\d$-Fourier expansion of $f^{\sharp}$ is congruent mod $p$ to
$$\sum_{(n,p)=1}\frac{a_n}{n} q^n-a_p \left(\sum a_m
q^{mp}\right) \frac{q'}{q^p}+\left( \sum a_m q^{mp^2}\right)
\left(\frac{q'}{q^p}\right)^p;$$ hence the Fourier expansion of
$f^{\sharp}$ is congruent mod $p$ to $\sum_{(n,p)=1}\frac{a_n}{n}
q^n$.}
\end{example}

\begin{example}
\cite{difmod} \textup{There exists a unique $\d$-modular form
$f^1\in H^0(X,L^{-\phi-1})$ with $\d$-Fourier expansion
$$
\sum_{n \geq 1} (-1)^{n-1}\frac{p^{n-1}}{n}\left(
\frac{q'}{q^p}\right)^n=\frac{q'}{q^p}+p\left( ...\right).$$
 Hence the Fourier expansion of $f^1$ is
$0$. By the way the above $\d$-Fourier expansion has the same shape
as a certain function playing a role in ``explicit" local class
field theory. Here is the rough idea for the construction of $f^1$. One considers the universal elliptic curve $E=\bigcup U_i\ra Spec\ M^{\infty}$, one considers sections $s_i:\widehat{U}_i\ra J^1(U_i)$ of the projection, one considers the differences
$$s_i-s_j:\widehat{U}_i \cap \widehat{U}_j \ra N^1\simeq \widehat{\mathbb G}_a,$$
and one considers the cohomology class $[s_i-s_j]\in H^1(\widehat{E},\cO)=H^1(E,\cO)$. Then $f^1\in M^{\infty}$ is defined as
the cup product of this class with the canonical generator of the $1$-forms; in fact $f^1$
is  the image of some element 
$$f^1(a_4,a_6,a'_4,a'_6)\in R[a_4,a_6,a'_4,a'_6,\Delta^{-1}]\h.$$
By the way it is a result of Hurlburt \cite{hurl} that 
$$f^1(a_4,a_6,a'_4,a'_6)\equiv E_{p-1}\frac{2a_4^p a'_6-3a_6^p a'_4}{\Delta^p}+f_0(a_4,a_6)\ \ \ \text{mod}\ \ \ p$$
where $f_0\in R[a_4,a_6,\Delta^{-1}]$ and we recall that $E_{p-1}$ is the Eisenstein form of weight $p-1$. 
(The polynomial $f_0$ is related to the Kronecker modular polynomial mod $p^2$.) This plus the  $\d$-Fourier expansion of $f^1$ should be viewed as an arithmetic analogue of the fact (due to Ramanujan) that
$$\frac{2a_4 da_6-3 a_6  da_4}{\Delta}\mapsto \frac{dq}{q}$$
via the map between the K\"{a}hler differentials induced by the Fourier expansion map $\bC[a_4,a_6,\Delta^{-1}]\ra \bC((q))$.} \end{example}

\begin{example}
 \cite{Barcau} \textup{There exists  unique $\d$-modular forms
$f^{\partial}\in H^0(X,L^{\phi-1})$ and $f_{\partial}\in H^0(X,L^{1-\phi})$ with $\d$-Fourier expansions $1$
(and hence Fourier expansions $1$). They satisfy $f^{\partial}\cdot f_{\partial}=1$.  The form $f^{\partial}$ is constructed by applying the ``top part" of  the canonical prolongation of the ``Serre operator" 
to the form $f^1$. More precisely $f^{\partial}$ is a constant times the image of
$$ (72 \phi(a_6)\frac{\partial}{\partial a'_4}-16\phi(a_4)^2\frac{\partial}{\partial a'_6}-p\cdot \phi(P))(f^1(a_4,a_6,a'_4,a'_6))$$
where $P\in R[a_4,a_6,\Delta^{-1},E_{p-1}^{-1}]\h$ is the Ramanujan form.}
\end{example}

\begin{theorem}
\label{uff}
\cite{Barcau, difmod, book}
The ring $R_{\d}(X,K^{-1})^{\sigma}$ is ``$\d$-generated" by $f^1$ and
$f^{\partial}$.
\end{theorem}

Note that $f^1$ and $f^{\partial}$ do not actually belong to
$R_{\d}(X,K^{-1})$  so  the  above theorem needs
some further explanation which we skip here; essentially, what happens is that $f^1$ and $f^{\partial}$ belong to a ring slightly bigger than
$R_{\d}(X,K^{-1})$ and they ``$\d$-generate" that ring.
We also note
the following   structure theorem
for the kernel and image of the $\d$-Fourier expansion map, in which the forms $f^1$ and $f^{\partial}$ play a key role:

\begin{theorem}
\label{vine} \cite{igusa}

1) The kernel of the Fourier expansion map  $M^{\infty}\ra R((q))\h$
is the $p$-adic closure of the smallest $\d$-stable ideal containing  $f^1$ and $f^{\partial}-1$.

2) The $p$-adic closure of the image of the Fourier expansion map
 $M^{\infty}\ra
R((q))\h$ equals Katz' ring $\bW$ of generalized $p$-adic modular
forms.
\end{theorem}

The proof of the Theorem above is rather indirect and heavily Galois-theoretic. 
Statement 2) in Theorem \ref{vine} says that all Katz' divided
congruences between classical modular forms can be obtained by
taking combinations of ``higher $p$-derivatives" of classical
modular forms. Statement 1) above is a lift to characteristic zero
of the Serre and Swinnerton-Dyer theorem about the kernel of the
classical Fourier expansion map for classical modular forms mod $p$.
Theorem \ref{uff} can also be viewed as a lift to characteristic zero of
results of Ihara \cite{Ihara1} about the Hasse invariant in
characteristic  $p$. These mod $p$ results do not lift to
characteristic zero in usual algebraic geometry but do lift, as we
see, to characteristic zero in $\d$-geometry.

We mention the following remarkable infinitesimal symmetry property; recall the classical Serre operator
$\partial:\cO(V)\ra \cO(V)$. Also consider the Euler derivation operator ${\mathcal D}:\cO(V)\ra \cO(V)$ given by multiplication by the degree on each graded component of $\cO(V)$. Finally let $P$ be the Ramanujan form (in the degree $2$ component of $\cO(V)$) and let $\theta:\cO(V)\ra \cO(V)$ be the derivation $\theta=\partial+P{\mathcal D}$.

\begin{theorem}
\cite{book}
The  operator $\theta$ is an infinitesimal symmetry
of
the $R$-module generated by $f^1,f^{\partial}$, and $f_{\partial}$ in $M^1=\cO^1(V)$. Also $\theta$ is a variational infinitesimal symmetry
of
the $R$-module generated by $f^{\partial}$ and $f_{\partial}$ in $M^1=\cO^1(V)$. 
\end{theorem}

Here is a calculation of differentials of the forms $f^1,f^{\partial},f_{\partial}$. 

\begin{theorem}
\cite{book}
Let $\omega,\alpha$ be the basis of $\Omega_{\cO(V)/R}$ dual to  $\theta, {\mathcal D}$. Then 
$$
\begin{array}{rcl}
\  & \  & \  \\
d (f^{\partial}) & = & f^{\partial} \cdot (\phi^* \alpha-\alpha)\\
\  & \  & \  \\
d(f_{\partial}) & = & -f_{\partial} \cdot (\phi^*\alpha- \alpha)\\
\  & \  & \  \\
d(f^1) & = & -f^1 \cdot (\phi^*\alpha+\alpha)-f_{\partial} \cdot \omega + f^{\partial}\cdot p^{-1}\phi^* \omega.
\end{array}
$$
In particular, for the eigenvalue $\mu=1$ we have
$$\int \{\frac{d(f^{\partial})}{f^{\partial}}\} dp=\int \{\frac{d(f_{\partial})}{f_{\partial}}\} dp=0.$$
\end{theorem}

By the way, the forms $f^{\sharp}$ and $f^1$ introduced above can be used
to prove some interesting purely diophantine results. For instance we have the following:

\begin{theorem}
\label{anti} \cite{local}
 Assume that $\Phi:X=X_1(N) \ra A$ is a modular
parametrization of an elliptic curve. Let $p$ be a sufficiently
large ``good" prime and let $Q \in X(R)$ be an ordinary point. Let
$S$ be the set of all rational primes that are inert in the
imaginary quadratic field attached to $Q$. Let $C$ be the
$S$-isogeny class of $Q$ in $X(R)$ (consisting of points
corresponding to isogenies of degrees only divisible by primes in
$S$). Then there exists a constant $c$ such that for any subgroup
$\Gamma \leq A(R)$ with $r:=rank(\Gamma)<\infty$ the set
$\Phi(C)\cap \Gamma$ is finite of cardinality at most $c p^{r}$.
\end{theorem}

Other results of the same type (e.g. involving Heegner points)  were
proved in \cite{local}. 
In particular an analogue of the above Theorem is true with $C$ replaced by the locus $CL$ of canonical lifts.
To have a rough idea about the arguments involved assume we want to prove that $\Phi(CL)\cap \Gamma$ is finite (and to bound the cardinality of this set) in case $\Gamma$ 
is the torsion group of $A(R)$. One considers the order $2$  $\d$-modular form $f^{\sharp}:X_1(N)(R)\ra R$ and one constructs, using $f^1$, a $\d$-function of order $1$,
$f^{\flat}:X(R)\ra R$,  on an open set $X \subset X_1(N)$ which vanishes exactly on $CL\cap X(R)$. Then any point $P$ in the intersection
 $X(R)\cap \Phi(CL)\cap \Gamma$ satisfies the system of ``differential equations of order $\leq 2$ in $1$ unknown" 
 $$\begin{cases}
 f^{\sharp}(P)=0\\
 f^{\flat}(P)=0
 \end{cases}$$
 One can show that this system is ``sufficiently non-degenerate" to allow the elimination of the ``derivatives" of the unknown; one is left
 with a differential equation $f^0(P)=0$ ``of order $0$" which has, then, only finitely many solutions (by Krasner's theorem). By that theorem one can also bound the number of solutions.
 
We end the discussion here by noting  that the main players in the theory above
enjoy a certain remarkable property which we call
$\d$-overconvergence. Morally this is a an overconvergence property (in the classical sense of Dwork, Monsky, Washnitzer) ``in the direction of the variables $x',x'',...,x^{(n)}$" (but not necessarily in the direction of $x$). We prove:

\begin{theorem} \cite{over} The
 $\d$-functions
$f^{\sharp}, f^1, f^{\partial}$  are  $\d$-overconvergent.
\end{theorem}

\subsection{$\d$-Hecke operators}
Next we discuss the   Hecke action on $\d$-modular forms.
For $(n,p)=1$ the Hecke operators $T_m(n)$ act naturally on
$\d$-series (i.e. series in $R((q))[q',...,q^{(r)}]\h$) by the usual formula that inserts roots
of unity of order prime to $p$ which are all in $R$. However no
naive definition of $T_m(p)$ seems to work. Instead we consider the
situation mod $p$ and make the following definition.
 Let $x$ be a $p$-tuple $x_1,...,x_p$ of indeterminates
and let $s$ be the $p$-tuple  $s_1,...,s_p$ of  fundamental
symmetric polynomials in $x$. An element $f\in
k[[q]][q',...,q^{(r)}]$ is called $\d$-$p$-symmetric mod $p$ if the
sum
$$f(x_1,...,x^{(r)}_1)+...+f(x_p,...,x^{(r)}_p)\in k[[x]][x',...,x^{(r)}]$$  is
the image of an element
$$f_{(p)}=f_{(p)}(s_1,...,s_p,...,s_1^{(r)},...,s_p^{(r)})\in
k[[s]][s',...,s^{(r)}].$$
 For $f$ that is $\d$-$p$-symmetric mod $p$ define
$$``pT_m(p)"f=f_{(p)}(0,...,0,q,...,0,...,0,q^{(r)}) +p^mf(q^p,...,\d^r(q^p))\in
k[[q]][q',...,q^{(r)}].$$ Eigenvectors of $``pT_m(p)"$ will be
automatically understood to be  $\d$-$p$-symmetric. Also let us say that a
series in $k[[q]][q',...,q^{(r)}]$ is primitive if the series in
$k[[q]]$ obtained by setting $q'=...=q^{(r)}=0$ is killed by the
classical $U$-operator.
Then one can give  a complete description (in terms of classical Hecke eigenforms mod $p$)
of $\d$-``eigenseries" mod $p$ of order $1$ which are $\d$-Fourier expansions of $\d$-modular forms of arbitrary order
and weight:

\begin{theorem}
\label{mooor} \cite{hecke} There is a $1-1$ correspondence between:

1) Series in $k[[q]]$ which are eigenvectors of all
$T_{m+2}(n),T_{m+2}(p)$, $(n,p)=1$, and which are Fourier expansions
of classical modular forms over $k$ of weight $\equiv m+2$ mod
$p-1$.

2) Primitive series in $k[[q]][q']$  which are eigenvectors of all
$nT_m(n),``pT_m(p)"$, $(n,p)=1$, and which are $\d$-Fourier
expansions of $\d$-modular forms of some order $\geq 0$ with weight
$w$, $deg(w)=m$.
\end{theorem}

Note that
the $\d$-Fourier expansion of the form $f^{\sharp}$ discussed in
Example \ref{cuc} is an example of series in 2) of  Theorem
\ref{mooor}. (Note that $f^{\sharp}$ has order $2$ although its
$\d$-Fourier expansion reduced mod $p$ has order $1$!) More
generally the series in 1) and 2) of Theorem \ref{mooor} are related
in an explicit way, similar to the way $f$ and $f^{\sharp}$ of
Example \ref{cuc} are related. The proof of  Theorem \ref{mooor}
involves a careful study of the action of $\d$-Hecke operators on $\d$-series plus the use of the canonical
prolongations of the Serre operator acting on $\d$-modular forms.

\subsection{$\d$-functions on finite flat schemes}
The $p$-jet spaces 
 of finite flat schemes over
$R$ seem to play a key role in many aspects of the theory. These $p$-jet spaces
 are neither finite nor flat in general and overall they seem quite pathological. 
 There are two remarkable classes of examples, however, where
 some order seems to be restored in the limit; these classes are finite flat $p$-group schemes that fit into $p$-divisible groups and finite length $p$-typical Witt rings. Recall that for any ring $A$ we write $\overline{A}=A/pA$. Then for connected $p$-divisible groups we have:

\begin{theorem}
\cite{pp} \label{mik}
 Let $\cF$ be a formal group law over $R$ in one variable $x$, assume ${\mathcal F}$ has  finite
height, and let $\cF[p^{\nu}]$ be the kernel of the multiplication by
$p^{\nu}$ viewed as a finite flat group scheme over $R$. Then
$$\lim_{\stackrel{\ra}{n}} \overline{\cO^n(\cF[p^{\nu}])}   \simeq
\frac{k[x,x',x'',...]}{(x^{p^{\nu}},(x')^{p^{\nu}},(x'')^{p^{\nu}},...)}$$
sending $x,\d x,\d^2 x,...$ into the classes of $x, x', x'',...$.
\end{theorem}

A similar result is obtained in \cite{pp} for the divisible group
$E[p^{\nu}]$ of an ordinary elliptic curve; some of the components
of $J^n(E[p^{\nu}])$ will be empty and exactly which ones are so is
dictated by the valuation of $q-1$ where $q$ is the Serre-Tate
parameter. The components that are non-empty (in particular the
identity component) behave in the same way as the formal groups
examined in Theorem \ref{mik} above.

In the same spirit one can compute $p$-jet spaces of Witt rings.
 Let us consider
 the ring $W_m(R)$ of $p$-typical Witt vectors of length $m+1$, $m\geq 1$, 
 and denote by $\Sigma_m=Spec\ W_m(R)$ its spectrum. 
 Set $v_i=(0,...,0,1,0,...,0)\in W_m(R)$, ($1$ preceded by $i$ zeroes, $i=1,...,m$), set $\pi=1-\d v_1\in \cO^1(\Sigma_m)$,   and let $\Omega_m=\{1,...,m\}$. The following is a description of the identity component of the limit of $p$-jet spaces  mod $p$:
 
 \begin{theorem}
 \label{uuu}
 \cite{pw}
 For $n \geq 2$ the image of $\pi^p$  in 
 $\overline{\cO^n(\Sigma_m)}$ is idempotent and we have an isomorphism
 $$\lim_{\stackrel{\ra}{n}}\overline{\cO^n(\Sigma_m)}_{\pi}\simeq \frac{k[x_i^{(r)};i\in \Omega_m;r \geq 0]}{(x_ix_j, (x_i^{(r)})^p;i,j\in \Omega_m, r\geq 1)}$$
 sending each $\overline{\d^r v_i}$ into the class of the variable $x_i^{(r)}$.
 \end{theorem}
 
 A similar description is obtained in \cite{pw} for the $p$-jet maps induced by the Witt comonad maps.
 We recall that the data consisting of $\cO^n(\Sigma_m)$ and the maps induced by the comonad maps should be viewed as 
 an arithmetic analogue of the Lie groupoid of the line.

\subsection{$\d$-Galois groups of $\d$-linear equations} Recall that for any solution $u\in GL_n(R)$ of a $\d$-linear equation $$\d u=\alpha \cdot u^{(p)}$$ (where $\alpha\in {\mathfrak gl}_n(R)$)  and for any $\d$-subring $\cO\subset R$ we defined the $\d$-Galois group $G_{u/\cO}\subset GL_n(\cO)$. We want to explain a result proved in \cite{adel2}.
Consider  the maximal torus $T\subset GL_n(R)$ of diagonal matrices, the Weyl group $W\subset GL_n(R)$ of permutation matrices,
the normalizer $N=WT=TW$ of $T$ in $GL_n(R)$, and the subgroup $N^{\d}$ of $N$ consisting of all elements of $N$ whose entries are in the monoid of constants $R^{\d}$.
We also use below the notation $K^a$ for the algebraic closure of the fraction field $K$ of $R$; the Zariski closed sets $Z$ of $GL_n(K^a)$ are then viewed as varieties over $K^a$.
The next result illustrates some ``generic" features of our $\d$-Galois groups; assertion 1) of the next theorem
shows that the $\d$-Galois
group is generically ``not too large". Assertions 2) and 3)
 show that the $\d$-Galois group are generically ``as large as possible". As we shall see presently, the meaning of the word {\it generic} is different in each of the $3$ situations: in situation 1) {\it generic} means {\it outside a Zariski closed set}; in situation 2) {\it generic} means {\it outside a thin set} (in the sense of diophantine geometry); in situation 3)  {\it generic} means {\it outside a set of the first category} (in the sense of Baire category).

\begin{theorem}
\label{food}\

1) There exists a Zariski closed subset $\Omega\subset GL_n(K^a)$ not containing $1$ such that for any $u\in GL_n(R) \backslash \Omega$ the following holds.  Let $\alpha=\d u \cdot (u^{(p)})^{-1}$ and let $\cO$ be a $\d$-subring of $R$ containing $\alpha$. Then
$G_{u/\cO}$ contains a normal subgroup of finite index which is diagonalizable over $K^a$.

2) Let $\cO=\bZ_{(p)}$. There exists a thin set $\Omega\subset {\mathbb Q}^{n^2}$ such that for any $\alpha\in \bZ^{n^2}\backslash \Omega$  there exists  a solution $u$ of the equation $\d u=\alpha u^{(p)}$ with the property that  $G_{u/\cO}$ is a finite group containing  the Weyl group $W$. 

3) There exists  a subset $\Omega$  of the first category in the metric space
$$X=\{u\in GL_n(R);u\equiv 1\ \ \text{mod}\ \ p\}$$
such that for any $u \in X\backslash \Omega$ the following holds.
 Let $\alpha=\d u\cdot (u^{(p)})^{-1}$. Then there exists a  $\d$-subring $\cO$ of $R$ containing $R^{\d}$ such that
$\alpha\in {\mathfrak gl}_n(\cO)$ and such that $G_{u/\cO}=N^{\d}$.
\end{theorem}

The groups $W$ and  $N^{\d}$ should be morally viewed as ``incarnations" of the groups  ``$GL_n({\mathbb F}_1)$" and ``$GL_n({\mathbb F}_1^a)$" where 
 ``${\mathbb F}_1$" and ``${\mathbb F}_1^a$" are the ``field with element" and ``its algebraic closure" respectively \cite{borger}. This suggests that the $\d$-Galois theory we are proposing here should be viewed as a Galois theory over ``${\mathbb F}_1$".

\subsection{$\d$-flows for the classical groups}
The main results in \cite{adel2} concern the existence of certain $\d$-flows on $\widehat{GL_n}$ that are adapted, in a certain precise sense, to the various classical subgroups $GL_n,SL_n,SO_n,Sp_n$. 

 Let $G=GL_n$ and
  let $H$  be a smooth closed subgroup scheme of $G$. We say that a $\d$-flow $\phi_G$ is 
 left  (respectively right) compatible with $H$ if $H$ is $\phi_{G}$-horizontal and
  the left (respectively right)  diagram below is commutative:
 $$
  \begin{array}{rcl}
\widehat{H}\times \widehat{G} & \ra & \widehat{G}\\
\phi_H\times \phi_G \downarrow & \  & \downarrow \phi_G\\
\widehat{H}\times \widehat{G} & \ra & \widehat{G}\end{array}\ \ \ \ 
\begin{array}{rcl}
\widehat{G}\times \widehat{H} & \ra & \widehat{G}\\
\phi_G\times \phi_H \downarrow & \  & \downarrow \phi_G\\
\widehat{G}\times \widehat{H} & \ra & \widehat{G}\end{array}
$$
where $\phi_H:\widehat{H}\ra \widehat{H}$ is induced by $\phi_G$ and the horizontal maps are given by multiplication.

  Next by an involution on $G$ we understand  a morphism of schemes $\dagger:G\ra G$ over $R$, $x\mapsto x^{\dagger}$, such that $x^{\dagger \dagger}=x$ and $(xy)^{\dagger}=y^{\dagger}x^{\dagger}$.
  By a quadratic map on $G$ we mean a morphism of schemes $\cH:G\ra G$ over $R$ such that there exists an involution $\dagger:G\ra G$ and an element $q\in G$ with the property that
  $\cH(x)=x^{\dagger}qx.$
  Given a quadratic map $\cH$ one can define a map
$ \cH_2:G\times G\ra G$ by the formula  $\cH_2(x,y)=\cH(x)x^{-1}y=x^{\dagger}qy$.
   One can also define a closed subgroup scheme
 $S$ of $G$ as the connected component $\cH^{-1}(q)^{\circ}$ of the group scheme $\cH^{-1}(q)$;
 we say that $S$ is defined by $\cH$. 

Let us fix now a lift of Frobenius  $\phi_{G,0}$ on $\widehat{G}$ such that $q$ is $\phi_{G,0}$-horizontal (i.e. the ideal of $q$ in $\cO(G)\h$ is sent by $\phi_{G,0}$ into itself). For a lift of Frobenius $\phi_{G}$ on $\widehat{G}$ we say that ${\mathcal H}$ is horizontal (respectively symmetric) with respect to  $\phi_G,\phi_{G,0}$ if 
 the left (respectively right) diagram below  is commutative:
  $$\begin{array}{rcl}
\widehat{G} & \stackrel{\phi_{G}}{\longrightarrow} & \widehat{G}\\
\cH  \downarrow & \  & \downarrow \cH \\
\widehat{G} & \stackrel{\phi_{G,0}}{\longrightarrow} & \widehat{G}\\
\end{array}\ \ \ \ 
\begin{array}{rcl}
\widehat{G} & \stackrel{\phi_{G,0} \times \phi_{G}}{\longrightarrow} & \widehat{G}\times \widehat{G}\\
\phi_{G} \times \phi_{G,0} \downarrow & \  & \downarrow \cH_2\\
 \widehat{G}\times \widehat{G} & \stackrel{\cH_2}{\longrightarrow} & \widehat{G}\end{array}
$$

Note that if ${\mathcal H}$ is horizontal with respect to $\phi_G,\phi_{G,0}$ and $q$ is $\phi_{G,0}$-horizontal  then the group $S$ defined by $\cH$ is $\phi_{G}$-horizontal; in particular there is an induced lift of Frobenius $\phi_S$ on $\widehat{S}$.
Also note that if we set  $\phi_{G,0}(x)=x^{(p)}$, viewing $\cH$ as a matrix $\cH(x)$ with entries in $R[x,\det(x)^{-1}]\h$, we have that horizontality of $\cH$ with respect to $\phi_{G},\phi_{G,0}$ is equivalent to the condition that
$\d_{G}\cH=0$, 
which can be interpreted as saying that ${\mathcal H}$ is a {\it prime integral} for our $\d$-flow $\phi_G$.

The basic split classical groups $GL_n,SO_n,Sp_n$ are defined by quadratic maps on $G=GL_n$ as follows. We start with
 $GL_n$ itself which is defined by  $\cH(x)=1$; in this case $x^{\dagger}=x^{-1}$, $q=1$. We call $\cH$ the canonical quadratic map defining $GL_n$. 
 We also recall that  $T\subset G$ is the maximal torus of diagonal matrices and
 $W\subset G$ is the Weyl subgroup of $GL_n$ of all permutation matrices.
Throughout our discussion we let $\phi_{G,0}(x)$ be the lift of Frobenius on $\widehat{GL_n}$ defined by $\phi_{G,0}(x):=x^{(p)}$; one can prove that this
 $\phi_{G,0}(x)$
  is the unique lift of Frobenius on $\widehat{G}$ that is left and right compatible with $T$ and $W$ and extends to a lift of Frobenius
 on $\widehat{{\mathfrak g}{\mathfrak l}_n}$ (where we view $\widehat{GL_n}$ as an open set of 
 $\widehat{{\mathfrak g}{\mathfrak l}_n}$).
On the other hand the groups $Sp_{2r}, SO_{2r}, SO_{2r+1}$ are defined by the quadratic map on $G=GL_n$ given by $\cH(x)=x^tqx$ where $q$ is equal to  
 $$
 \left(\begin{array}{cl} 0 & 1_r\\-1_r & 0\end{array}\right),\ \ 
\left( 
\begin{array}{ll} 0 & 1_r\\1_r & 0\end{array}\right),\ \ 
\left( \begin{array}{lll} 1 & 0 & 0\\
0 & 0 & 1_r\\
0 & 1_r & 0\end{array}\right),\ \ 
$$
 $n=2r, 2r, 2r+1$ respectively, $x^{\dagger}=x^t$ is the transpose, and $1_r$ is the $r\times r$ identity matrix. We call this $\cH$ the canonical quadratic map defining $Sp_{2r}, SO_{2r}, SO_{2r+1}$ respectively. All these groups are smooth over $R$.

   \begin{theorem}\label{laugh}
 Let $S$ be any  of the groups $GL_n, SO_n, Sp_n$ and let
 $\cH$ be the canonical quadratic map on $G=GL_n$ defining $S$.
 Then the following hold.
 
 1) (Symmetry and horizontality.) There exists a  unique lift of Frobenius $\phi_{G}$ on $\widehat{G}$ such that ${\mathcal H}$  is horizontal and symmetric with respect to $\phi_{G},\phi_{G,0}$.
 
  2) (Compatibility with torus and Weyl group.) $\phi_{G}$   is right compatible with $T$ and $W$; also $\phi_{G}$ is  left compatible with $T_S:=T\cap S$ and $W_S:=W\cap S$. In particular if $l\d:GL_n\ra {\mathfrak g}{\mathfrak l}_n$ is the arithmetic logarithmic derivative associated to $\phi_{G}$ then for all $a\in  T_S W_S$  and $b\in GL_n$ (alternatively for all $a\in GL_n$ and $b\in TW$) we have
  $$
\label{cris}l\d(ab)=a \star_{\d} l\d (b) +_{\d} l\d(a).$$
  
    3) (Compatibility with root groups.) If $\chi$ is a root of $S$ (which is not a shortest root of $SO_n$ with $n$ odd) then the
    corresponding root group $U_{\chi}\simeq {\mathbb G}_a$ is $\phi_{GL_n}$-horizontal.

 \end{theorem} 
 
Note that a similar result can be proved for $SL_n$; the involution $\dagger$ lives, in this case, on a cover of $GL_n$ rather than on $GL_n$ itself. Note also that the exception in assertion 3 of the Theorem (occurring in case $\chi$ is a shortest root of $SO_n$ with $n$ odd) is a curious phenomenon which deserves further understanding.

\section{Problems}

\begin{problem}
\label{p1}
 Study the arithmetic jet spaces $J^n(X)$ of curves $X$
(and more general varieties) with bad reduction.\end{problem}

 This
could be applied, in particular,  to tackle  Mazur's question
\cite{Mazur} about bounding the torsion points on curves unifromly
in terms of the genus; in other words replacing $C(g,p)$ by $C(g)$
in Theorem \ref{zuzu}.
 Our proof in \cite{pjets} is based on
the study of the arithmetic jet space of $J^1(X)$ at a prime $p$ of
good reduction. A study of the arithmetic jet space of curves at
primes of bad reduction might lead to dropping the dependence of
$C(g,p)$ on $p$. Evidence that arithmetic jet spaces can be handled
in the case of bad reduction comes in particular from the recent
paper \cite{over}.

\begin{problem}
\label{p2}
 Study the $\d$-modular forms  that vanish on
arithmetically interesting Zariski dense subsets of Shimura
varieties  (such as CM loci or individual non-CM isogeny classes).
Compute $\d$-invariants of higher dimensional correspondences.
\end{problem}

 This could be applied to extend results in
\cite{local}, e.g. Theorem \ref{anti} above. A deeper study of
differential modular forms may allow one, for instance, to replace
$S$-isogeny class with the full isogeny class. The arguments might
then be extended to higher dimensional contexts and to the global
field rather than the local field situation. That such a deeper
study is possible is shown by papers like \cite{igusa}, for
instance. For the higher dimensional case the theory in
\cite{siegel} might have to be further developed to match the one
dimensional theory in \cite{difmod, book}. In a related direction
one might attempt to use the methods in \cite{local} to tackle
Pink's conjectures in \cite{Pink}. In \cite{local} it was shown that
behind finiteness theorems in diophantine geometry one can have
reciprocity maps that are somehow inherited from $\d$-geometry (and
that provide effective bounds); a similar picture might hold for
(cases of) Pink's conjecture. The first case to look at for such
reciprocity maps would be in the case of the intersection between a
multisection $X$ of an abelian (or semiabelian) scheme $G\ra S$ over
a curve $S$ with the set of torsion points lying in special (CM or
otherwise) fibers; more general situations, in which torsion points
are replaced by division points of a group generated by finitely
many sections, can be considered. Results of Andr\'{e}, Ribet, and
Bertrand are pertinent to this question.

\begin{problem}
\label{p3} Compare the $\d$-geometric approach to quotient spaces
with the approach via non-commutative geometry. \end{problem}

The quotients $X/\sigma$ for the correspondences appearing in
Theorem \ref{thm2} do not exist, of course, in usual algebraic
geometry. As Theorem \ref{thm2} shows these quotients exist,
however,  and are interesting in $\d$-geometry. Remarkably such
quotients also exist and are interesting in non-commutative geometry
\cite{Marcolli}. More precisely  the $3$ cases (spherical, flat,
hyperbolic) of Theorem \ref{thm2}
 are closely related to the following $3$ classes of examples studied in
  non-commutative geometry:

   1) (spherical) $\frac{\bP^1(\bR)}{SL_2(\bZ)}$, non-commutative
   boundary of the classical modular curve;

2) (flat) $\frac{S^1}{\langle e^{2 \pi i
   \tau}\rangle}$ ($\theta \in \bR\backslash \bQ$):
non-commutative elliptic curves;

3) (hyperbolic) Non-commutative space $Sh^{nc}$ containing the
classical Shimura variety $Sh$ ($2$-dimensional analogue of
Bost-Connes systems).

 It would be interesting to understand
     why these $3$ classes appear in both contexts ($\d$-geometry and non-commutative geometry);
     also one would like to see
      whether there is a connection, in the case of these $3$ classes, between the $2$ contexts.
      
      Note that non-commutative geometry can also tackle the dynamics of rational functions  that are not necessarily post-critically finite of Euler characteristic zero.  It is conceivable that some post-critically finite polynomials of non-zero Euler characteristic
      possess $\d$-invariants for some particular  primes (with respect to the anticanonical bundle or other bundles).
      A good start would be to investigate the $\d$-invariants of $\sigma(x)=x^2-1$. Another good start would be to investigate $\d$-invariants of post-critically finite polynomials with Euler characteristic zero that are congruent modulo special primes to post-critically finite polynomials with Euler characteristic zero.

\begin{problem}
\label{p4} Study the de Rham cohomology  of arithmetic jet spaces.
Find arithmetic analogues of K\"{a}hler differentials $\Omega$ and
${\mathcal D}$-modules. Find an object that is to ${\mathcal D}$
what $\cO^1$ is to $Sym(\Omega)$.
\end{problem}

 The study of de Rham cohomology of arithmetic jet spaces
  was started in \cite{forms} where it is
shown that the de Rham cohomology of  $J^n(X)$ carries information
about the arithmetic of $X$. The de Rham computations in
\cite{forms} are probably shadows of more general phenomena which
deserve being understood. Also the de Rham setting could be replaced
by an overconvergent one; overconvergence is known to give an
improved picture of the de Rham story and, on the other hand, as already mentioned, it was
proved in \cite{over} that most of the remarkable $\d$-functions
occurring in the theory possess a remarkable overconvergence
property in the ``$\d$-variables" called $\d$-{\it overconvergence}.
Finally one is tempted to try to relate the de Rham cohomology of
arithmetic jet spaces $J^n(X)$ to the de Rham-Witt complex of $X$ in
characteristic $p$ and in mixed characteristic. Note further that
since the arithmetic jet space $J^1(X)$ is an analogue of the
(physical) tangent bundle $T(X)$ of $X$ it follows that the sheaf
$\cO^1$ is an arithmetic analogue of the sheaf
$Sym(\Omega_{X/R})$, symmetric algebra on the K\"{a}hler
differentials. But there is no obvious arithmetic analogue of the
sheaf $\Omega_{X/R}$ itself. Also there is no obvious arithmetic
analogue of the sheaf ${\mathcal D}_X$ of differential operators on
$X$ and of ${\mathcal D}_X$-modules. The absence of immediate
analogues of $\Omega$ and ${\mathcal D}$ is of course related to the
intrinsic non-linearity of $p$-derivations. It would be interesting
to search for such analogues. It is on the other hand conceivable
that there is a sheaf in the arithmetic theory that is to ${\mathcal
D}$ what $\cO^1$ is to $Sym(\Omega)$. Recall that the associated
graded algebra of ${\mathcal D}$ is canonically isomorphic to the
algebra of functions on the (physical) cotangent bundle $\cO(T^*(X))$
(and not on the tangent bundle); this looks like a discrepancy but
actually the arithmetic jet space $J^1(X)$ has a sort of intrinsic
self-duality (cf. \cite{book}) that is missing in the classical
algebro-geometric case where the tangent bundle $T(X)$ and the
cotangent bundle $T^*(X)$ and not naturally dual (unless, say, a
symplectic structure is given).

The $\d$-overconvergence property mentioned in Problem \ref{p4} may
hold the key to:

\begin{problem}
\label{p5}
 Define and study  the/a  maximal space of $\d$-modular forms on
which Atkin's $U$ operator can be defined. \end{problem}

Indeed the Hecke operators $T(n)$ with $p\not|n$ are defined on
$\d$-modular forms and have a rich theory in this context
\cite{difmod}. In contrast to this $T(p)$ and hence $U$ are still
mysterious in the theory of $\d$-modular forms. An  step in
understanding $U$ was taken in \cite{hecke} where the theory mod $p$
for series of order $1$  was given a rather definitive treatment.
However the theory in characteristic zero seems elusive at this
point. There are two paths towards such a theory so far. One path is
via $\d$-{\it symmetry} \cite{dcc,hecke}; this is a characteristic
$0$ analogue of the concept of $\d$-$p$-symmetry mod $p$ discussed above.
  Another path is via $\d$-overconvergence
\cite{over}. The two paths seem to lead into different directions
and this discrepancy needs to be better understood. Assuming that a
good theory of $U$ is achieved, this might lead to a Hida-like
theory of families of differential modular forms, including  Galois
representations attached to such forms. It is conceivable that
families in this context are not power series but Witt vectors.
Part of the quest for a $U$ theory of $\d$-modular forms is to seek a $\d$-analogue
of Eisentein series. It is conceivable that the rings $\cO^n(X_1(N))$ contain functions that do not vanish at the cusps and are eigenvectors of the Hecke operators; such functions could be viewed as ``$\d$-Eisenstein" forms of weight zero.

\begin{problem}
\label{p6}
  Interpret information contained in the arithmetic jet spaces $J^n(X)$
  as an arithmetic Kodaira-Spencer ``class" of
  $X$.\end{problem}

   Indeed some of these arithmetic Kodaira-Spencer classes (e.g
  in the case of elliptic curves or, more generally, abelian schemes)
  were studied in \cite{difmod,book} and lead to interesting $\d$-modular forms.
For general schemes (e.g for curves of higher genus) these classes were explored in Dupuy's thesis \cite{Taylor}. They
are non-abelian cohomology classes with values in the sheaf of
automorphisms of $p$-adic affine spaces $\widehat{\mathbb A}^d$ (in
the case of curves $d=1$). These classes arise from comparing the
local trivializations of arithmetic jet spaces. In this more general case these
classes  may hold the key to a ``deformation theory over the field
with one element''.  On the other hand Dupuy proved in \cite{Taylor} that if $X$ is a smooth projective curve of genus $\geq 2$ over $R$ then $J^1(X)$ is a torsor for some line bundle over $X$; this is rather surprising in view of the high non-linearity of the theory.
One should say that the line bundle in question is still mysterious  and deserves further investigation.

\begin{problem}
\label{p7} Further develop the partial differential theory in
\cite{laplace,pde,pdemod}.\end{problem}

Indeed in spite of the extensive work done in
\cite{laplace,pde,pdemod} the arithmetic {\it partial} differential
theory is still in its infancy. The elliptic case of that theory
\cite{laplace} (which, we recall, involves operators $\delta_{p_1}$
and $\delta_{p_2}$ corresponding to two primes $p_1$ and $p_2$) is
directly related to the study of the de Rham cohomology of
arithmetic jet spaces \cite{forms}; indeed one of the main results
in \cite{forms} shows that the arithmetic Laplacians in
\cite{laplace} are formal primitives (both $p_1$-adically and
$p_2$-adically) of global $1$-forms on the arithmetic jet spaces
(these forms being not exact, although formally exact, and hence
closed). By the way analogues of these results in \cite{forms}
probably exist in the case of modular curves; in the one prime case
a beginning of such a study was undertaken in \cite{book}, where
some of the main $\d$-modular forms of the theory were shown to
satisfy some remarkable systems of Pfaff equations. The
hyperbolic/parabolic case of the theory \cite{pde,pdemod}  (which, we
recall, involves a $p$-derivation $\d_p$ with respect to a prime $p$
and a usual derivation operator $\delta_q$) could be further developed
 as follows. One could start by ``specializing" the
variable $q$ in $\delta_q$ to elements $\pi$ in {\it arbitrarily ramified}
extensions of $\bZ_p$. This might push the theory in the ``arbitrarily ramified
direction" which would be extremely desirable for
arithmetic-geometric applications. Indeed our ordinary arithmetic
differential theory is, at present, a non-ramified (or at most ``boundedly ramified") theory. A further idea along these lines would be to use the solutions of
 the arithmetic partial differential equations in \cite{pde,pdemod} to
 let points ``flow" on varieties defined over number fields.
Some of the solutions in \cite{pde,pdemod} have interesting
arithmetic features (some look like hybrids between quantum
exponentials and Artin-Hasse exponentials, for instance) so the
``flows" defined by them might have arithmetic consequences. The
challenge is to find (if at all possible) ``special values" of these
solutions that are algebraic. One should also mention that the
arithmetic hyperbolic and parabolic equations in \cite{pde,pdemod}
have, in special cases, well defined ``indices" that seem to carry
arithmetic information; the challenge would be to make the index
machinery work in general situations and to study the variation of
indices in families.

\begin{problem}
\label{p8} Find an arithmetic analogue of Sato hyperfunction
solutions of both ``ordinary" and ``partial" arithmetic differential
equations. \end{problem}

Indeed Sato's hyperfunctions, in their simplest incarnation, are
pairs $(f(x),g(x))$ of functions on the unit disk (corresponding to
the distribution $f(x)-g(x^{-1})$) modulo $(c,c)$, $c$ a constant.
The derivative of a pair is then
$(\frac{df}{dx}(x),-x^2\frac{dg}{dx}(x))$.) One could then try to
consider, in the arithmetic case, pairs $(P,Q)$ of points of
algebraic groups with values in $\d$-rings  modulo an appropriate
equivalence relation and with an appropriate analogue of
differentiation with respect to $p$; this framework could be the
correct one for ``non-analytic"  solutions of the equations in
\cite{laplace,pde,pdemod}.

\begin{problem}
\label{p9} Construct $p$-adic measures from $\d$-modular forms.
\end{problem}

Indeed  one of the main ideas in Katz's approach to $p$-adic
interpolation \cite{Katzgen} was to lift some of the remarkable
$\bZ_p$-valued ($p$-adic) measures of the theory to ${\mathbb
W}$-valued  measures where ${\mathbb W}$ is the ring of
(generalized) $p$-adic modular forms. One can hope that some of
these ${\mathbb W}$-valued measures of Katz can be further lifted to
measures with values in the $p$-adic completion of the ring  of
$\d$-modular functions $M^{\infty}$. Indeed recall from Theorem
\ref{vine} that there is a canonical homomorphism $M^{\infty}\ra
{\mathbb W}$ whose image is $p$-adically dense, hence the ``lifting"
problem makes sense. These lifted $\widehat{M^{\infty}}$-valued
measures could then be evaluated at various elliptic curves defined
over $\d$-rings  to obtain new $\bZ_p$-valued measures (and hence
new $p$-adic interpolation results)  in the same way in which Katz
evaluated his measures at special elliptic curves. Another related
idea would be to interpret the solutions in $\bZ_p[[q]]$ of the
arithmetic partial differential equations in Problem \ref{p7} above
as measures (via Iwasawa's representation of measures as power
series). There is a discrepancy in this approach in that the
derivation of interest in Iwasawa's theory is $(1+q)\frac{d}{dq}$
whereas the derivation of interest in Problem \ref{p7} is
$q\frac{d}{dq}$; nevertheless one should pursue this idea and
understand the discrepancy.

\begin{problem}
\label{p10} Find arithmetic analogues of classical theorems in the
theory of differential algebraic groups and further develop the $\d$-Galois theory in \cite{adel3}.
\end{problem}

 Indeed the theory of groups defined by (usual) differential equations
(``differential algebraic groups'') is by now a classical subject:
it goes back to Lie and Cartan and underwent a new development, from
a rather new angle, through the work of Cassidy and Kolchin
\cite{Cassidy,Kolchin}. It is tempting to seek an arithmetic
analogue of this theory: one would like to understand, for instance,
the structure of all subgroups of $GL_n(R)$ that are defined by
arithmetic differential equations. The paper \cite{buca} proves an
arithmetic analogue of Cassidy's theorem about Zariski dense
subgroups of simple algebraic groups over differential fields; in
\cite{buca} the case of Zariski dense mod $p$ groups is considered.
But Zariski dense groups such as $GL_n(\bZ_p)$ are definitely
extremely interesting (and lead to interesting Galois theoretic
results such as in \cite{book}, Chapter 5). So a generalization of
\cite{buca} to the case of Zariski dense (rather than Zariski
dense mod $p$) groups, together with a generalization of the Galois
theoretic results in \cite{book}, would be very desirable. For instance 
it would be interesting to classify all $\delta$-subgroups of the multiplicative group
${\mathbb G}_m(R)=R^{\times}$ (or more generally of $GL_2(R)$) and find the invariants of such groups acting on 
$\widehat{R\{x\}_{(p)}}$ (where $R\{x\}=R[x,x',x'',...]$). Also remark that, as in the case of usual derivations,
there are interesting (``unexpected") homomorphisms in $\d$-geometry between $GL_1={\mathbb G}_m$ and
$GL_2$, for instance $${\mathbb G}_m(R)\rightarrow GL_2(R),\ \ \ \ a \mapsto \left( \begin{array}{rcl} a & a\psi_*(a)\\ 0 & a \end{array}\right)$$
where $\psi_*$ is a $\d$-character, i.e. $\psi \in {\mathcal X}^n({\mathbb G}_m)$.
Cf. \cite{herras} for interesting developments into this subject.
The  main open problem in the $\d$-Galois theory of $GL_n$ \cite{adel3} seems at this point  to decide if the $\d$-Galois groups always contain a subgroup of the diagonal matrices as a subgroup of finite index. Other problems are:  to establish a Galois correspondence;  to understand the relation (already hinted at in \cite{adel3}) between $\d$-Galois groups and Galois problems arising from the dynamics on ${\mathbb P}^n$;  to generalize the theory by replacing $GL_n$ with an arbitrary reductive group.

\begin{problem}
\label{p11} Compose some of the basic $\d$-functions of the theory
in \cite{book} (e.g. the $\d$-modular forms on Shimura curves) with
$p$-adic uniformization maps (e.g. with Drinfeld's uniformization
map of Shimura curves). \end{problem}

Indeed this might shed a new light (coming from the analytic world)
on $\d$-geometry. These composed maps would belong to a ``$\d$-rigid
geometry" which has yet to be developed in case these maps are
interesting enough to require it. By the way it is not at all clear
that the functor that attaches to a formal $p$-adic scheme its
arithmetic jet space can be prolonged to a functor in the rigid
category. The problem (which seems to boil down to some quite
non-trivial combinatorially flavored calculation) is to show that
the arithmetic jet space functor sends blow-ups with centers in the
closed fiber into (some version of) blow-ups.

In the spirit of the last comments on Problem \ref{p11} one can ask:

\begin{problem}
\label{p12} Study the morphisms $J^n(X)\ra J^n(Y)$  induced by
non-\'{e}tale finite flat covers of smooth schemes $X \ra
Y$.\end{problem}

 Indeed note that if $X \ra Y$ is finite and \'{e}tale then it is well known that
 $J^n(X) \ra J^n(Y)$ are finite and \'{e}tale; indeed $J^n(X)\simeq J^n(Y)\times_Y X$.
 But note that if $X \ra Y$ is only finite and flat then $J^n(X)\ra J^n(Y)$ is
 neither finite nor flat in general.
  One of the simplest
examples which need to be investigated (some partial results are
available \cite{pp}, cf. Theorem \ref{mik}) is that of the covers
$[p^{\nu}]:G\ra G$ of smooth group schemes (or formal groups) given
by multiplication by $p^{\nu}$. The geometry of the induced
endomorphisms of the arithmetic jet spaces is highly complex and
mysterious. Understanding it might be, in particular, another path
towards introducing/understanding the Atkin operator $U$ on
$\d$-modular forms referred to in Problem \ref{p5}. An obviously
closely related problem is to understand the $p$-jet spaces of
$p$-divisible groups; cf. Theorem \ref{mik}. Yet another example of
interest is the study of $J^n(X)\ra J^n(X/\Gamma)$ for $\Gamma$ a
finite group acting on a smooth $X$; even the case $X=Y^n$ with $Y$
a curve and $\Sigma=S_n$ the symmetric group acting naturally is
still completely mysterious. This latter case is related to the
concept of $\d$-symmetry mentioned in Problem \ref{p5} and appeared
in an essential manner in the paper \cite{dcc}: the failure of
invariants to commute with formation of jet spaces (in the ramified
case) is directly responsible for the existence (and indeed
abundance) of $\d$-functions on smooth projective curves which do
not arise from $\d$-characters of the Jacobian.

\begin{problem}
\label{p13} Understand which analytic functions $X(\bZ_p)\ra \bZ_p$
for smooth schemes $X/\bZ_p$ are induced by $\d$-functions.
\end{problem}

Indeed  it was proved in \cite{BRS} that  a function $f:\bZ_p
\rightarrow \bZ_p$ is analytic if, and only if, there exists $m$
such that $f$ can be represented as $f(x)=F(x, \d x, \ldots , \d^m
x)$, where $F\in \bZ_p[x_0,x_1,...,x_m]\h$ is a restricted power series with
$\bZ_p$-coefficients in $m+1$ variables.  This can be viewed as a
``differential interpolation result": indeed $f(x)$ is given by a
finite family of (unrelated) power series $F_i(x)$ convergent on
disjoint balls $B_i$ that cover $\bZ_p$ and the result says that one
can find a single power series $F(x,\d x,...,\d^m x)$ that equals
$F_i(x)$ on $B_i$ for each $i$. One can ask for a generalization of
this by asking which analytic functions $f:X(\bZ_p)\ra \bZ_p$
defined on the $\bZ_p$-points of a smooth scheme $X/\bZ_p$ come from
a $\d$-function $\tilde{f} \in \cO^n(X)$ (i.e. $f=\tilde{f}_*$); of
course such a $\tilde{f}$ cannot be, in general,  unique. If
$X=\bA^1$ is the affine line then the result in \cite{BRS} says that
any $f$ comes from some $\tilde{f}$. This is probably still the case
if $X$ is any affine smooth scheme. On the other hand this fails if $X=\bP^1$ is
the projective line simply because there are no non-constant
$\d$-functions in $\cO^n(\bP^1)$ \cite{pjets} but, of course, there are plenty
of non-constant analytic functions $\bP^1(\bZ_p)\ra \bZ_p$. There
should be a collection of  cohomological
obstructions to lifting $f$ to some $\tilde{f}$ that should reflect
the global geometry of $X$. This seems to us a rather fundamental
question in understanding the relation between $p$-adic analytic
geometry and $\d$-geometry.

In the light of Theorem \ref{mik} one can ask:

\begin{problem}
\label{p14} Let $\alpha_0,...,\alpha_{\nu}$ be elements in the
algebraic closure of the fraction field of $R$ which are integral
over $R$. Let $X_i=Spec\ R[\alpha_i]$ be viewed as a closed
subscheme of the line $\bA^1$ over $R$. Compute/understand the
arithmetic jet spaces $J^n(\cup_{i=0}^{\nu} X_i)$.
\end{problem}

This is a problem ``about the interaction"  of algebraic numbers.
Indeed the case $\nu=0$ is clear; for instance if $R[\alpha]$ is
totally ramified over $R$ (and $\neq R$) then the arithmetic jet
spaces are empty: $J^n(Spec\ R[\alpha])=\emptyset$ for $n \geq 1$.
However, for $\nu \geq 1$, an interesting new phenomenon occurs.
Indeed if $\alpha_i=\zeta_{p^i}$ (primitive $p^i$-th root of unity)
then although $J^n(Spec\ R[\zeta_{p^i}])=\emptyset$ for $n\geq 1$
and $i\geq 1$ we have that $J^n(\cup_{i=0}^{\nu} Spec\
R[\zeta_{p^i}])=J^n(\mu_{p^{\nu}})$ is non-empty and indeed
extremely interesting; cf. Theorem \ref{mik}.

\begin{problem}
\label{p15}  Find
arithmetic analogues of  Hamiltonian systems and of algebraically completely integrable systems. Find arithmetic analogues
of the formal pseudo-differential calculus.
\end{problem}

This problem  is motivated by the link (due to Fuchs and Manin
\cite{Manin}) between the Painlev\'{e} VI equation (which has a Hamiltonian structure) and the Manin map of an
elliptic fibration. Painlev\'{e} VI possesses an arithmetic analogue whose study was begun in \cite{BYM}; this study is in its infancy and deserves further attention.
 More generally there  is  an intriguing
possibility that other physically relevant differential equations
(especially arising in Hamiltonian contexts, especially in the
completely integrable situation, both finite and infinite dimensional) have arithmetic analogues carrying
arithmetic significance. Finally it is conceivable that a meaningful arithmetic analogue
of the formal pseudo-differential calculus in one variable \cite{olver}. p. 318 can be developed;
in other words one should be able to bring into the picture  negative powers
of the $p$-derivation $\d:\cO^n(X)\ra \cO^{n+1}(X)$ in the same way in which formal pseudo-differential calculus brings into the picture the negative powers of the total derivative operator acting on functions on jet spaces.

\begin{problem}
\label{p16} Find an arithmetic analogue of the differential
groupoids/Lie pseudogroups  in the work of Lie, Cartan,  Malgrange \textup{[39]}.
\end{problem}

Indeed recall that
one can view
 the $p$-jet spaces $J^n(\Sigma_m)$ of $\Sigma_m=Spec\ W_m(R)$ and the jet maps induced by the comonad maps
 as an arithmetic analogue of the Lie groupoid $J^n(\bR\times \bR/\bR)^*$ of the line; cf. \cite{pw} and  Theorem \ref{uuu} above for results on this. It would be interesting to investigate subobjects of the system $J^n(\Sigma_m)$ that play the role of analogues of differential sub-groupoids and to find arithmetic analogues of the differential invariants of diffeomorphism groups acting on natural bundles arising from frame bundles. 
 A theory of $J^n(\Sigma_m)$, suitably extended to other rings and to arithmetic analogues of poly-vector fields, could also be interpreted as an arithmetic analogue of the theory of the deRham-Witt complex \cite{he} because it would replace the usual K\"{a}hler differentials $\Omega$ by constructions involving the operators $\d_p$.

  \begin{problem}
  Compute $\cO^n({\mathbb T}(N,\kappa))$ where ${\mathbb T}(N,\kappa)$ is the Hecke ${\mathbb Z}$-algebra
  attached to cusp forms on $\Gamma_1(N)$ of level $\kappa$. 
  \end{problem}
  
  This could lead to a way of associating differential modular forms $f^{\sharp}$ to classical newforms $f$ of weight $\neq 2$. (Cf. Theorem \ref{mooor}.) This might also lead to a link between $\d$-geometry and Galois representations. The problem  may be related to the study of double coset sets of $GL_2$ (or more generally $GL_n$) and a link of this to the study in 
  \cite{adel1,adel2,adel3} is  plausible.
  
  \medskip
  
 We end by
stating two of the most puzzling concrete open problems of the theory.

\begin{problem}
\label{p17} Compute $\cO^n(X)$ for an elliptic curve $X$.
\end{problem}

Note that $\cO^1(X)$ and $\cX^n(X)$ have been computed (Theorem \ref{thm1}) and one
expects $\cX^n(X)$ to ``generate" $\cO^n(X)$. In particular assume $X$ is
not a canonical lift and $\psi$ is a basis of $\cX^2(X)$; is it true
that $\cO^n(X)=R[\psi,\d \psi,...,\d^{n-1}\psi]\h$? The analogue of
this in differential algebra is true; cf. \cite{hermann}.

\begin{problem}
\label{p18} Compute $\cX^n(G)$ for $G$ an extension  of an elliptic
curve by $\bG_m$.  \end{problem}

One expects that this module depends in an arithmetically
interesting way on the class of the extension. The problem is
directly related to that of understanding the cohomology of
arithmetic jet spaces.

\end{document}